\documentclass{article}
\usepackage{amssymb}

\usepackage{graphicx}
\usepackage{amsmath}

\newtheorem{theorem}{Theorem}

\newtheorem{corollary}[theorem]{Corollary}

\newtheorem{definition}[theorem]{Definition}

\newtheorem{lemma}[theorem]{Lemma}

\newtheorem{proposition}[theorem]{Proposition}

\input{tcilatex}

\begin{document}

\title{Mirror symmetry without corrections}
\author{Naichung Conan Leung\thanks{%
This paper is partially supported by a NSF grant, DMS-9803616.} \\
School of Mathematics, \\
University of Minnesota, \\
Minneapolis, MN 55455, U.S.A..}
\date{}
\maketitle

\begin{abstract}
We give geometric explanations and proofs of various mirror symmetry
conjectures for $T^{n}$-invariant Calabi-Yau manifolds when instanton
corrections are absent. This uses fiberwise Fourier transformation together
with base Legendre transformation.

We discuss mirror transformations of

(i) moduli spaces of complex structures and complexified symplectic
structures, $H^{p,q}$'s, Yukawa couplings;

(ii) $\mathbf{sl}\left( 2\right) \times \mathbf{sl}\left( 2\right) $ actions;

(iii) holomorphic and symplectic automorphisms and

(iv) A- and B-connections, supersymmetric A- and B-cycles, correlation
functions.

We also study (ii) for $T^{n}$-invariant hyperkahler manifolds.
\end{abstract}

\tableofcontents

\newpage

Mirror symmetry conjecture predicts that there is a transformation from
complex (resp. symplectic) geometry of one Calabi-Yau manifold $M$ to
symplectic (resp. complex) geometry of another Calabi-Yau manifold $W$ of
the same dimension. Such pairs of manifolds are called mirror manifolds.
This transformation should also has the inversion property, namely if we
take the transformation twice, we recover the original geometry.

It is expected that such transformation exists for Calabi-Yau manifolds near
a \textit{large complex structure limit} point. Such point in the moduli
space should correspond to the existence of a semi-flat Calabi-Yau metric,
possibly highly singular.

To understand \textit{why and how }these two different kinds of geometry got
interchanged between mirror manifolds, we study the $T^{n}$-invariant case
in details. The importance of the $T^{n}$-invariant (or more generally
semi-flat) case is first brought up by Strominger, Yau and Zaslow in their
foundational paper \cite{SYZ} which explains mirror symmetry from a
physical/geometric viewpoint. This is now called the SYZ mirror conjecture.
The $T^{n}$-invariant case is then studied by Hitchin in \cite{H1}, Yau,
Zaslow and the author in \cite{LYZ} and it is also an important part of this
paper. The main advantage here is the absence of holomorphic disks, the
so-called instantons.

We start with an affine manifold $D$, which we assume to be a domain in $%
\mathbb{R}^{n}$ in this introduction. Let $\phi $ be an elliptic solution to
the real Monge-Amper\'{e} equation on $D$: 
\begin{eqnarray*}
\det \nabla ^{2}\phi  &=&1, \\
\nabla ^{2}\phi  &>&0.
\end{eqnarray*}
Then it determines two noncompact Calabi-Yau manifolds, $TD$ and $T^{\ast }D$%
. Notice that $T^{\ast }D$ carries a canonical symplectic structure and $TD$
carries a canonical complex structure because $D$ is affine. We can also
compactify the fiber directions by quotienting $TD$ and $T^{\ast }D$ with a
lattice $\Lambda $ in $\mathbb{R}^{n}$ and its dual lattice $\Lambda ^{\ast }
$ in $\mathbb{R}^{n\ast }$ respectively and obtain mirror manifolds $M$ and $%
W$. The natural fibrations of $M$ and $W$ over $D$ are both special
Lagrangian fibrations.

The mirror transform from $M$ to $W,$ and vice versa, is basically (i) the
Fourier transformation on fibers of $M\rightarrow D$ together with (ii) the
Legendre transformation on the base $D$. The Calabi-Yau manifold $W$ can
also be identified as the moduli space of flat $U\left( 1\right) $
connections on special Lagrangian tori on $M$ with its $L^{2}$ metric. We
are going to explain how the mirror transformation exchanges complex
geometry and symplectic geometry between $M$ and $W$:

(1) The identification between moduli spaces of complex structures on $M$
and complexified symplectic structures on $W$, moreover this map is both
holomorphic and isometric;

(2) The identification of $H^{p,q}\left( M\right) $ and $H^{n-p,q}\left(
W\right) $;

(3) The mirror transformation of certain A-cycles in $M$ to B-cycles in $W$.
We also identify their moduli spaces and correlation functions (this is
partly borrowed from \cite{LYZ}). In fact the simplest case here is the
classical Blaschke connection and its conjugate connection, they got
interchanged by Legendre transformation;

(4) There is an $\mathbf{sl}\left( 2\right) $ action on the cohomology of $M$
induced from variation of Hodge structures. Together with the $\mathbf{sl}%
\left( 2\right) $ action from the hard Lefschetz theorem, we obtain an $%
\mathbf{sl}\left( 2\right) \times \mathbf{sl}\left( 2\right) $ action on the
cohomology of $M$. Under mirror transformation from $M$ to $W$, these two $%
\mathbf{sl}\left( 2\right) $ actions interchange their roles;

(5) Transformations of holomorphic automorphisms of $M$ to symplectic
automorphisms of $W$, in fact its preserves a naturally defined two tensor
on $W$, not just the symplectic two form.

\bigskip

In the last section we study $T^{n}$-invariant hyperk\"{a}hler manifolds.
That is when the holonomy group of $M$ is inside $Sp\left( n/2\right)
\subset SU\left( n\right) $. Cohomology of a hyperk\"{a}hler manifold admits
a natural $\mathbf{so}\left( 4,1\right) $ action. In the $T^{n}$-invariant
case, we show that our $\mathbf{sl}\left( 2\right) \times \mathbf{sl}\left(
2\right) =\mathbf{so}\left( 3,1\right) $ action on cohomology is part of
this hyperk\"{a}hler $\mathbf{so}\left( 4,1\right) $ action.

In \cite{KS} Kontsevich and Soibelman also study mirror symmetry for these $%
T^{n}$-invariant Calabi-Yau manifolds, their emphasis is however very
different from ours.

Acknowledgments: The author thanks Richard Thomas, Xiaowei Wang, Shing-Tung
Yau and Eric Zaslow for many helpful and valuable discussions. The author
also thanks Mark Gross for pointing out an earlier mistake on B-fields and
other comments. The paper is prepared when the author visited the Natural
Center of Theoretical Science, Tsing-Hua University, Taiwan in the summer of
2000. The author thanks the center for providing an excellent research
environment and support. This project is also partially supported by a NSF
grant, DMS-9803616.

\section{$T^{n}$-invariant Calabi-Yau and their mirrors}

A Calabi-Yau manifold $M$ of real dimension $2n$ is a Riemannian manifold
with $SU\left( n\right) $ holonomy, or equivalently a K\"{a}hler manifold
with zero Ricci curvature. We can reduce this condition to a complex
Monge-Amper\'{e} equation provided that $M$ is compact. Yau proved that this
equation is always solvable as long as $c_{1}\left( M\right) =0$, vanishing
of the first Chern class of $M$.

Even though it is easy to construct Calabi-Yau manifolds, it is extremely
difficult to write down their Ricci flat metrics. When Calabi-Yau manifolds
have $T^{n}$ symmetry, we can study translation invariant solution to the
complex Monge-Amper\'{e} equation and reduces the problem to the solution of
a real Monge-Amper\'{e} equation.

These $T^{n}$-invariant Calabi-Yau manifolds is a natural class of semi-flat
Calabi-Yau manifolds. Recall that a Calabi-Yau manifold is called semi-flat
if it admits a fibration by flat Lagrangian tori. Such manifolds are
introduced into mirror symmetry in \cite{SYZ} and then further studied in 
\cite{H1}, \cite{Gr} and \cite{LYZ}.

\bigskip 

\bigskip 

\textbf{The real Monge-Amper\'{e} equation}

First we consider the dimension reduction of the complex Monge-Amper\'{e}
equation to the real Monge-Amper\'{e} equation. The resulting Ricci flat
metric would be a $T^{n}$-invariant Calabi-Yau metric: Let $M$ be a tubular
domain in $\mathbb{C}^{n}$ with complex coordinates $z^{j}=x^{j}+iy^{j}$, 
\begin{equation*}
M=D\times i\mathbb{R}^{n}\subset \mathbb{C}^{n},
\end{equation*}
where $D$ is a convex domain in $\mathbb{R}^{n}$. The holomorphic volume
form on $M$ is given by 
\begin{equation*}
\Omega _{M}=dz^{1}\wedge dz^{2}\wedge \cdots \wedge dz^{n}\text{.}
\end{equation*}
Let $\omega _{M}$ be the K\"{a}hler form of $M$, then the complex
Monge-Amper\'{e} equation for the Ricci flat metric is the following: 
\begin{equation*}
\Omega _{M}\bar{\Omega}_{M}=C\omega _{M}^{n}\text{.}
\end{equation*}

We assume that the K\"{a}hler potential $\phi $ of the K\"{a}hler form $%
\omega _{M}=i\partial \bar{\partial}\phi $ is invariant under translations
along imaginary directions. That is, 
\begin{equation*}
\phi \left( x^{j},y^{j}\right) =\phi \left( x^{j}\right)
\end{equation*}
is a function of the $x^{j}$'s only. In this case the complex
Monge-Amper\'{e} equation becomes the real Monge-Amper\'{e} equation. Cheng
and Yau \cite{CY} proved that there is a unique elliptic solution $\phi
\left( x\right) $ to the corresponding boundary value problem 
\begin{eqnarray*}
\medskip \det \left( \frac{\partial ^{2}\phi }{\partial x^{j}\partial x^{k}}%
\right) &=&C,\, \\
\,\,\phi |_{\partial D} &=&0.
\end{eqnarray*}
Ellipticity of a solution $\phi $ is equivalent to the convexity of $\phi $,
i.e. 
\begin{equation*}
\left( \frac{\partial ^{2}\phi }{\partial x^{j}\partial x^{k}}\right) >0%
\text{.}
\end{equation*}

\bigskip

We can compactify imaginary directions by taking a quotient of $i\mathbb{R}%
^{n}$ by a lattice $i\Lambda $. That is we replace the original $M$ by $%
M=D\times iT$ where $T$ is the torus $\mathbb{R}^{n}/\Lambda $ and the above
K\"{a}hler structure $\omega _{M}$ descends to $D\times iT$. If we write 
\begin{equation*}
\phi _{jk}=\frac{\partial ^{2}\phi }{\partial x^{j}\partial x^{k}},
\end{equation*}
then the Riemannian metric on $M$ is 
\begin{equation*}
g_{M}=\Sigma \phi _{jk}\left( dx^{j}\otimes dx^{k}+dy^{j}\otimes
dy^{k}\right)
\end{equation*}
and the symplectic form $\omega _{M}$ is 
\begin{equation*}
\omega _{M}=\frac{i}{2}\Sigma \phi _{jk}dz^{j}\wedge d\bar{z}^{k}.
\end{equation*}

Notice that $\omega _{M}$ can also be expressed as 
\begin{equation*}
\omega _{M}=\Sigma \phi _{jk}dx^{j}\wedge dy^{k}
\end{equation*}
because of $\phi _{jk}=\phi _{kj}$. The closedness of $\omega _{M}$ follows
from $\phi _{ijk}=\phi _{kji}=\partial _{i}\partial _{j}\partial _{k}\phi $.

Remark: It is easy to see that $D\times i\Lambda \subset TD$ is a special
Lagrangian submanifold (for the definition of a special Lagrangian, readers
can refer to later part of this section.)

\bigskip

\textbf{Affine manifolds and complexifications}

Notice that the real Monge-Amper\'{e} equation 
\begin{equation*}
\det \left( \frac{\partial ^{2}\phi }{\partial x^{j}\partial x^{k}}\right)
=const,
\end{equation*}
is invariant under any affine transformation 
\begin{equation*}
\left( x^{j}\right) \rightarrow \left( \bar{x}^{j}\right) =\left(
A_{k}^{j}x^{k}+B^{j}\right) .
\end{equation*}
This is because 
\begin{equation*}
\frac{\partial ^{2}\phi }{\partial x^{j}\partial x^{k}}=A_{j}^{l}A_{k}^{m}%
\frac{\partial ^{2}\phi }{\partial \bar{x}^{l}\partial \bar{x}^{m}},
\end{equation*}
and 
\begin{equation*}
\det \left( \frac{\partial ^{2}\phi }{\partial x^{j}\partial x^{k}}\right)
=\det \left( A\right) ^{2}\det \left( \frac{\partial ^{2}\phi }{\partial 
\bar{x}^{l}\partial \bar{x}^{m}}\right) .
\end{equation*}

The natural spaces to study such equation are affine manifolds. A manifold $%
D $ is called an affine manifold if there exists local charts such that
transition functions are all affine transformations as above. Over $D$,
there is a natural real line bundle whose transition functions are given by $%
\det A$. We denote it by $\mathbb{R\rightarrow }L\rightarrow D$. Now if $%
\phi \left( x\right) $ is a solution to the above equation with $const=1$ on
the coordinate chart with local coordinates $x^{j}$'s. Under the affine
coordinate change $\bar{x}^{j}=A_{k}^{j}x^{k}+B^{j}$, the function $\bar{\phi%
}=\left( \det A\right) ^{2}\phi $ satisfies 
\begin{equation*}
\det \left( \frac{\partial ^{2}\bar{\phi}}{\partial \bar{x}^{j}\partial \bar{%
x}^{k}}\right) =1.
\end{equation*}
Therefore on a general affine manifold $D$, a solution to the real
Monge-Amper\'{e} equation (with $const=1$) should be considered as a section
of $L^{\otimes 2}$.

\bigskip

It is not difficult to see that the tangent bundle of an affine manifold is
naturally an affine complex manifold: If we write a tangent vector of $D$ as 
$\Sigma y^{j}\frac{\partial }{\partial x^{j}}$ locally, then $%
z^{j}=x^{j}+iy^{j}$'s are local holomorphic coordinates of $TD$. The
transition function for $TD$ becomes $\left( z^{j}\right) \rightarrow \left(
A_{k}^{j}z^{k}+B^{j}\right) $, hence $TD$ is an affine complex manifold.

We want to patch the $T^{n}$-invariant Ricci flat metric on each coordinate
chart of $TD$ to the whole space and thus obtaining a $T^{n}$-invariant
Calabi-Yau manifold $M=TD$ (or $TD/\Lambda $). To do this we need to assume
that $\det A=1$ for all transition functions, such $D$ is called a special
affine manifold. Then 
\begin{eqnarray*}
g_{M} &=&\Sigma \phi _{jk}\left( dx^{j}\otimes dx^{k}+dy^{j}\otimes
dy^{k}\right) \\
\omega _{M} &=&\Sigma \phi _{jk}\left( x\right) dx^{j}\wedge dy^{k}=\frac{i}{%
2}\Sigma \phi _{jk}dz^{j}\wedge d\bar{z}^{k}.
\end{eqnarray*}
are well-defined K\"{a}hler metric and K\"{a}hler form over the affine
complex manifold $M$, which has a fibration over the real affine manifold $D$%
. Moreover 
\begin{equation*}
g_{D}=\Sigma \phi _{jk}\left( x\right) dx^{j}\otimes dx^{k}
\end{equation*}
defines a Riemannian metric on $D$ of Hessian type.

\bigskip

\textbf{Legendre transformation}

All our following discussions work for $D$ being a special orthogonal affine
manifold. For simplicity we assume that $D$ is simply a convex domain in $%
\mathbb{R}^{n}$ and $M=TD=D\times i\mathbb{R}^{n}$.

It is well-known that one can produce another solution to the real
Monge-Amper\'{e} equation from any given one via the so-called Legendre
transformation: We consider a change of coordinates $x_{k}=x_{k}\left(
x^{j}\right) $ given by 
\begin{equation*}
\frac{\partial x_{k}}{\partial x^{j}}=\phi _{jk},
\end{equation*}
thanks to the convexity of $\phi $. Then we have 
\begin{equation*}
\frac{\partial x^{j}}{\partial x_{k}}=\phi ^{jk},
\end{equation*}
where 
\begin{equation*}
\left( \phi ^{jk}\right) =\left( \phi _{jk}\right) ^{-1}\text{.}
\end{equation*}
Since $\phi ^{jk}=\phi ^{kj}$, locally there is a function $\psi \left(
x_{k}\right) $ on the dual vector space $\mathbb{R}^{n\ast }$ such that 
\begin{equation*}
x^{j}\left( x_{k}\right) =\frac{\partial \psi \left( x_{k}\right) }{\partial
x_{j}}\text{.}
\end{equation*}
Therefore, 
\begin{equation*}
\phi ^{jk}=\frac{\partial ^{2}\psi }{\partial x_{j}\partial x_{k}}\text{.}
\end{equation*}
This function $\psi \left( x_{k}\right) $ is called the Legendre
transformation of the function $\phi \left( x^{j}\right) $. It is obvious
that the convexity of $\phi $ and $\psi $ are equivalent to each other.
Moreover 
\begin{equation*}
\det \left( \frac{\partial ^{2}\phi }{\partial x^{j}\partial x^{k}}\right)
=C,
\end{equation*}
is equivalent to

\begin{equation*}
\det \left( \frac{\partial ^{2}\psi }{\partial x_{j}\partial x_{k}}\right)
=C^{-1}.
\end{equation*}
Furthermore the Legendre transformation has the inversion property, namely
the transformation of $\psi $ is $\phi $ again.

\bigskip

\textbf{Dual tori fibration - fiberwise Fourier transformation}

This construction works for any $T^{n}$-invariant K\"{a}hler manifold $M$,
not necessary a Calabi-Yau manifold. On $M=D\times iT$ there is a natural
torus fibration structure given by the projection to the first factor, 
\begin{eqnarray*}
M &\rightarrow &D, \\
\left( x^{j},y^{j}\right) &\rightarrow &\left( x^{j}\right) .
\end{eqnarray*}
Instead of performing the Legendre transformation to the base of this
fibration, we are going to replace the fiber torus $T=\mathbb{R}^{n}/\Lambda 
$ by the dual torus $T^{\ast }=\mathbb{R}^{n\ast }/\Lambda ^{\ast },$ where $%
\Lambda ^{\ast }=\left\{ v\in \mathbb{R}^{n\ast }:v\left( u\right) \in 
\mathbb{Z}\text{ for any }u\in \Lambda \right\} $ is the dual lattice to $%
\Lambda $.

In dimension one, taking the dual torus is just replacing a circle of radius 
$R$ to one with radius $1/R$. In general, if $y^{j}$'s are the coordinates
for $T$ and $y_{j}$'s their dual coordinates. Then a flat metric on $T$ is
given by $\Sigma \phi _{jk}dy^{j}\otimes dy^{k}$ for some constant positive
definite symmetric tensor $\phi _{jk}$. As usual we write 
\begin{equation*}
\left( \phi ^{jk}\right) =\left( \phi _{jk}\right) ^{-1},
\end{equation*}
then $\Sigma \phi ^{jk}dy_{j}\otimes dy_{k}$ is the dual flat metric on $%
T^{\ast }$.

\bigskip

Now we write $W=D\times iT^{\ast }$, the fiberwise dual torus fibration to $%
M=D\times iT$. Since the metric $g_{M}$ on $M$ is $T^{n}$-invariant, its
restriction to each torus $\left\{ x\right\} \times iT$ is the flat metric $%
\Sigma \phi _{jk}\left( x\right) dy^{j}\otimes dy^{k}$. The dual metric on
the dual torus $\left\{ x\right\} \times iT^{\ast }$ is $\phi ^{jk}\left(
x\right) dy_{j}\otimes dy_{k}$. So the natural metric on $W$ is given by 
\begin{equation*}
g_{W}=\Sigma \phi _{jk}dx^{j}\otimes dx^{k}+\phi ^{jk}dy_{j}\otimes dy_{k}.
\end{equation*}
If we view $T^{\ast }$ as the moduli space of flat $U\left( 1\right) $
connections on $T$, then it is not difficult to check that the
Weil-Petersson $L^{2}$ metric on $T^{\ast }$ is also $\phi
^{jk}dy_{j}\otimes dy_{k}$.

If we ignore the lattice structure for the moment, then $M$ is the tangent
bundle $TD$ of an affine manifold and 
\begin{equation*}
W=T^{\ast }D,
\end{equation*}
moreover, $g_{W}$ is just the induced Riemannian metric on the cotangent
bundle from the Riemannian metric $g_{D}=\Sigma \phi _{jk}dx^{j}\otimes
dx^{k}$ on $D$.

\bigskip

Even though $T^{\ast }D$ does not have a natural complex structure like $TD$%
, it does carry a natural symplectic structure: 
\begin{equation*}
\omega _{W}=\Sigma dx^{j}\wedge dy_{j},
\end{equation*}
which is well-known and plays a fundamental role in symplectic geometry. $%
\omega _{W}$ and $g_{W}$ together determine an almost complex structure $%
J_{W}$ on $W$ as follow, 
\begin{equation*}
\omega _{W}\left( X,Y\right) =g_{W}\left( J_{W}X,Y\right) .
\end{equation*}
In fact this almost complex structure is integrable and the holomorphic
coordinates are given by $z_{j}=x_{j}+iy_{j}$'s where $x_{j}\left( x\right) $
is determined by the Legendre transformation $\frac{\partial x_{j}}{\partial
x^{k}}=\phi _{jk}$ as before. In terms of this coordinate system, we can
rewrite $g_{W}$ and $\omega _{W}$ as follows 
\begin{eqnarray*}
g_{W} &=&\Sigma \phi ^{jk}\left( dx_{j}\otimes dx_{k}+dy_{j}\otimes
dy_{k}\right) \\
\omega _{W} &=&\frac{i}{2}\Sigma \phi ^{jk}dz_{j}\wedge d\bar{z}_{k}.
\end{eqnarray*}

\bigskip

Suppose that $g_{M}$ is a Calabi-Yau metric on $M$, namely $\phi \left(
x^{j}\right) $ satisfies the real Monge-Amper\'{e} equation, then $\psi
\left( x_{j}\right) $ also satisfies the real Monge-Amper\'{e} equation
because of 
\begin{equation*}
\phi ^{jk}=\frac{\partial ^{2}\psi }{\partial x_{j}\partial x_{k}}.
\end{equation*}
Therefore the metric $g_{W}$ on $W$ is again a $T^{n}$-invariant Calabi-Yau
metric.

\bigskip

We call this combination of the Fourier transform on fibers and the Legendre
transform on the base of a $T^{n}$-invariant K\"{a}hler manifold the \textit{%
mirror transformation.}

The similarities between $g_{M},\omega _{M}$ and $g_{W},\omega _{W}$ are
obvious. In particular, the mirror transformation has the inversion
property, namely the transform of $W$ is $M$ again.

\bigskip

Here is an important observation: On the tangent bundle $M=TD$, suppose we
vary its symplectic structure while keeping its natural complex structure
fixed. We would be looking at a family of solutions to the real
Monge-Amper\'{e} equation. On the $W=T^{\ast }D$ side, the corresponding
symplectic structure is unchanged, namely $\omega _{W}=\Sigma dx^{j}\wedge
dy_{j}$. But the complex structures on $W$ varies because the complex
coordinates on $W$ are given by $dz_{j}=\phi _{jk}dx^{k}+idy_{j}$ which
depends on particular solutions of the real Monge-Amper\'{e} equation.

By the earlier remark about the symmetry between $M$ and $W$, changing the
complex structures on $M$ is also equivalent to changing the symplectic
structures on $W$. To make this precise, we need to consider complexified
symplectic structures by adding B-fields as we will explain later.

In fact the complex geometry and symplectic geometry of $M$ and $W$ are
indeed interchangeable! String theory predicts that such phenomenon should
hold for a vast class of pairs of Calabi-Yau manifolds. This is the famous
Mirror Symmetry Conjecture.

General Calabi-Yau manifolds do not admit $T^{n}$-invariant metrics,
therefore we want to understand the process of constructing $W$ from $M$ via
a geometric way. To do this we need to introduce A- and B-cycles.

\bigskip

\textbf{Supersymmetric A- and B-cycles}

It was first argued by Strominger, Yau and Zaslow \cite{SYZ} from string
theory considerations that the mirror manifold $W$ should be identified as
the moduli space of special Lagrangian tori together with flat $U\left(
1\right) $ connections on them. These objects are called supersymmetric
A-cycles (see for example \cite{MMMS}, \cite{L1}). Let us recall the
definitions of A-cycles and B-cycles (we also include the B-field in these
definitions, see the next section for discussions on B-fields).

\begin{definition}
Let $M$ be a Calabi-Yau manifold of dimension $n$ with complexified K\"{a}%
hler form $\omega ^{\mathbb{C}}=\omega +i\beta $ and holomorphic volume form 
$\Omega $. We called a pair $\left( C,E\right) $ a \textit{supersymmetric
A-cycle} (or simply \textit{A-cycle)}, if (i) $C$ is a special Lagrangian
submanifold of $M$, namely $C$ is a real submanifold of dimension $n$ with 
\begin{equation*}
\omega |_{C}=0,
\end{equation*}
and 
\begin{equation*}
\func{Im}e^{i\theta }\Omega |_{C}=0,
\end{equation*}
for some constant angle $\theta $ which is called the \textit{phase angle}.

(ii) $E$ is a unitary vector bundle on $C$ whose curvature tensor $F$
satisfies the deformed flat condition, 
\begin{equation*}
\beta |_{C}+F=0\text{.}
\end{equation*}
\end{definition}

Note that the Lagrangian condition and deformed flat equation can be
combined into one complex equation on $C$: 
\begin{equation*}
\omega ^{\mathbb{C}}+F=0.
\end{equation*}

\begin{definition}
Let $M$ be a K\"{a}hler manifold with complexified K\"{a}hler form $\omega ^{%
\mathbb{C}}$, we called a pair $\left( C,E\right) $ a \textit{supersymmetric
B-cycle} (or simply \textit{B-cycle)}, if $C$ is a complex submanifold of $M$
of dimension $m$, $E$ is a holomorphic vector bundle on $C$ with a Hermitian
metric whose curvature tensor $F$ satisfies the following \textit{deformed
Hermitian-Yang-Mills equations }on $C$: 
\begin{equation*}
\func{Im}e^{i\theta }\left( \omega ^{\mathbb{C}}+F\right) ^{m}=0,
\end{equation*}
for some constant angle $\theta $ which is called the \textit{phase angle}.
\end{definition}

\bigskip

Remark: The following table gives a quick comparison of these two kinds of
supersymmetric cycles (see \cite{L1} for details).

\begin{equation*}
\left. 
\begin{array}{c}
\,\text{A-cycles} \\ 
E\rightarrow C\subset M \\ 
\\ 
\medskip \omega ^{\mathbb{C}}+F=0 \\ 
\\ 
\func{Im}e^{i\theta }\Omega =0
\end{array}
\right| 
\begin{array}{c}
\text{B-cycles} \\ 
E\rightarrow C\subset M \\ 
\\ 
\iota _{\Lambda ^{k-1}T_{M}}\Omega +F^{2,0}=0 \\ 
\\ 
\func{Im}e^{i\theta }\left( \omega ^{\mathbb{C}}+F\right) ^{m}=0.
\end{array}
\end{equation*}

\bigskip

\textbf{Constructing the mirror manifold}

Now we consider the moduli space of A-cycles $\left( C,E\right) $ on $M$
with $C$ a torus and the rank of $E$ equals one. In \cite{SYZ} SYZ
conjecture that $W$ is the mirror manifold of $M$. The L$^{2}$ metric on
this moduli space is expected to coincide with the Calabi-Yau metric on $W$
after suitable \textit{corrections} which comes from contributions from
holomorphic disks in $M$ whose boundaries lie on these A-cycles, these are
called instantons.

When $M$ is a $T^{n}$-invariant Calabi-Yau manifold with fibration $\pi
:M\rightarrow D$ as before. Then each fiber of $\pi $ is indeed a special
Lagrangian torus and $D$ is their moduli space. Each fiber together with the
restricted metric is a flat torus. Its dual torus can be naturally
identified with the moduli space of flat $U\left( 1\right) $ connections on
it. Therefore the space $W$, obtained by replacing each fiber torus in $M$
by its dual, can be naturally identified as the moduli space of A-cycles in
this case.

Furthermore the $L^{2}$ metric on this moduli space coincides with the dual
metric $g_{W}$ up to a constant multiple. Physically this is because of the
absence of instanton in this case. We have the following simple result.

\begin{theorem}
Under the natural identification of $W$ with the moduli space of flat $%
U\left( 1\right) $ connections on special Lagrangian tori in $M$, the metric 
$g_{W}$ equals the L$^{2}$ metric on the moduli space multiply with the
volume of the fiber.
\end{theorem}

Proof: Recall that $D$ is the moduli space of special Lagrangian tori. Let $%
\frac{\partial }{\partial x^{j}}$ be a tangent vector at a point in $D$, say
the origin. This corresponds to a harmonic one form on the central fiber $%
C\subset M$. This harmonic one form on $C$ is $\Sigma \phi _{jk}\left(
0\right) dy^{k}$. Now the moduli space L$^{2}$ inner product of $\frac{%
\partial }{\partial x^{j}}$ and $\frac{\partial }{\partial x^{l}}$ equals 
\begin{eqnarray*}
&\ll &\frac{\partial }{\partial x^{j}},\frac{\partial }{\partial x^{l}}\gg \\
&=&\int_{C}\left\langle \phi _{jk}dy^{k},\phi _{lm}dy^{m}\right\rangle dv_{C}
\\
&=&\int_{C}\phi _{jk}\phi _{lm}\phi ^{km}dv_{C} \\
&=&\phi _{jl}\left( 0\right) vol\left( C\right) \text{.}
\end{eqnarray*}
On the other hand 
\begin{equation*}
g_{W}\left( \frac{\partial }{\partial x^{j}},\frac{\partial }{\partial x^{l}}%
\right) =\phi _{jl}\left( 0\right) .
\end{equation*}
Similarly we can identify metrics along fiber directions of $\pi
:W\rightarrow D$. By definition the L$^{2}$ metric has no mixed terms
involving both the base and fiber directions. Hence we have the theorem. $%
\square $

\bigskip

\textbf{Shrinking the torus fibers}

Now we fix the symplectic form on $M$ as $\omega _{M}=\Sigma \phi
_{jk}dx^{j}\wedge dy^{k}$ and vary the complex structures. Instead of using
holomorphic coordinates $z^{j}=x^{j}+iy^{j}$'s, we define the new complex
structure on $M$ using the following holomorphic coordinates, 
\begin{equation*}
z_{t}^{j}=\frac{1}{t}x^{j}+iy^{j},
\end{equation*}
for any $t\in \mathbb{R}_{>0}$. The corresponding Calabi-Yau metric becomes 
\begin{equation*}
g_{t}=\Sigma \phi _{jk}\left( \frac{1}{t}dx^{j}\otimes dx^{k}+tdy^{j}\otimes
dy^{k}\right) .
\end{equation*}
The same fibration $\pi :M\rightarrow D$ is a special Lagrangian fibration
for each $t$. Moreover the volume form on $M$ is independent of $t$, namely $%
dv_{M}=\omega _{M}^{n}/n!$. As $t$ goes to zero, the size of the fibers
shrinks to zero while the base gets infinitely large.

If we rescale the metric to $tg_{t}$, then the diameter of $M$ stays bound
and $\left( M,tg_{t}\right) $'s converge in the Gromov-Hausdorff sense to
the real $n$ dimension manifold $D$ with the metric $g_{D}=\Sigma \phi
_{jk}dx^{j}dx^{k}$ as $t$ approaches zero. It is expected that similar
behaviors hold true for Calabi-Yau metrics near the large complex structure
limit, as least over a large portion of $M$. This prediction is verified by
Gross and Wilson when $M$ is a K3 surface \cite{GW}.

$\bigskip $

$\mathbf{B}$\textbf{-fields}

The purpose of introducing $B$-fields is to complexify the space of
symplectic structures on $M$, the conjectural mirror object to the space of
complex structures on $W$ which is naturally a complex space. Readers could
skip this part for the first time.

The usual definition of a B-field $\beta $ is a harmonic two form of type $%
\left( 1,1\right) $ on $M$, i.e. $\beta \in \Omega ^{1,1}\left( M,\mathbb{R}/%
\mathbb{Z}\right) $ with $d\beta =0$ and $d^{\ast }\beta =0$. These are
equivalent to the following conditions, 
\begin{eqnarray*}
d\beta  &=&0 \\
\beta \wedge \omega _{M}^{n-1} &=&c^{\prime }\omega _{M}^{n}\text{.}
\end{eqnarray*}
It is shown by Gross \cite{Gr} that if we consider a closed form $\beta $,
then the modified Legendre transformation, as we will describe later, does
preserve the Calabi-Yau condition on the $W$ side. However the harmonicity
will be loss. To remedy this problem, we need to deform the harmonic
equation. There are two natural way to do this, depending on whether we
prefer the complex polarization or the real polarization. We will first
discuss the one using the real polarization, namely the special Lagrangian
fibration. When we are in the large complex and K\"{a}hler structure limit,
the complex conjugation is the same as the real involution which sends the
fiber directions to its negative, namely $dx^{j}\rightarrow dx^{j}$ and $%
dy^{j}\rightarrow -dy^{j}$ in our previous coordinates. In general the
complex conjugation and the real involution are different. We denote the
holomorphic volume form of $M$ under this real involution (resp. complex
conjugation) by $\widehat{\Omega }$ (resp. $\bar{\Omega}$). The author
thanks Gross for pointing out an earlier mistake about $\widehat{\Omega }$.
Just like the distinction between K\"{a}hler metrics and those satisfying
the Monge-Amper\'{e} equations, namely Calabi-Yau metrics, we need the
following definitions.

\begin{definition}
Let $M$ be a Calabi-Yau manifold with holomorphic volume form $\Omega $.
Suppose that $\omega $ is a K\"{a}hler form on $M$ and $\beta $ is a closed
real two from on $M$ of type $\left( 1,1\right) $. Then $\omega +i\beta $ is
called a complexified Calabi-Yau K\"{a}hler form on $M$ and we denote it $%
\omega ^{\mathbb{C}}$ if $\left( \omega ^{\mathbb{C}}\right) ^{n}$ is a
nonzero constant multiple of $i^{n}\Omega \wedge \widehat{\Omega }$. 
\begin{equation*}
\left( \omega ^{\mathbb{C}}\right) ^{n}=ci^{n}\Omega \wedge \widehat{\Omega }%
\text{.}
\end{equation*}
We call this the complexified complex Monge-Amper\'{e} equation.
\end{definition}

\bigskip

\bigskip

The second definition of a B-field is to use $\bar{\Omega}$ and require that
the Calabi-Yau manifold $M$ satisfies 
\begin{eqnarray*}
\omega ^{n} &=&ci^{n}\Omega \bar{\Omega} \\
\func{Im}e^{i\theta }\left( \omega +i\beta \right) ^{n} &=&0 \\
\func{Im}e^{i\phi }\Omega &=&0\text{ on the zero section.}
\end{eqnarray*}

If we expand the second equation near the large K\"{a}hler structure limit,
namely we replace $\omega $ by a large multiple of $\omega $, or
equivalently we replace $\beta $ by a small multiple of it, we have 
\begin{equation*}
\left( \omega +i\varepsilon \beta \right) ^{n}=\omega ^{n}+i\varepsilon
n\beta \omega ^{n-1}+O\left( \varepsilon ^{2}\right) \text{.}
\end{equation*}
So if we linearize this equation, by deleting terms of order $\varepsilon
^{2}$ or higher. Then it becomes 
\begin{equation*}
\beta \omega ^{n-1}=c^{\prime }\omega ^{n}\text{.}
\end{equation*}
That is $\beta $ is a harmonic real two form. This approximation is in fact
the usual convention for a B-field.

\bigskip

\textbf{Including }$B$\textbf{-fields in the $T^{n}$-invariant case}

We first consider the case when $\omega ^{\mathbb{C}}=\omega +i\beta $
satisfies the complexified Monge-Amper\'{e} equation.

We suppose $\pi :M\rightarrow D$ is a $T^{n}$-invariant Calabi-Yau manifold
as before and $\omega _{M}=\Sigma \phi _{jk}\left( x\right) dx^{j}dy^{k}$ is
a $T^{n}$-invariant K\"{a}hler form on it. As usual we will include a $B$%
-field on $M$ which is invariant along fiber directions of $\pi $, namely $%
\beta _{M}=i\partial \bar{\partial}\eta \left( x\right) $. It is easy to see
that 
\begin{equation*}
\beta _{M}=\frac{i}{2}\Sigma \eta _{jk}\left( x\right) dz^{j}\wedge d\bar{z}%
^{k}=\Sigma \eta _{jk}\left( x\right) dx^{j}\wedge dy^{k}
\end{equation*}
with 
\begin{equation*}
\eta _{jk}=\eta _{kj}=\frac{\partial ^{2}\eta }{\partial x^{j}\partial x^{k}}%
\text{.}
\end{equation*}
Then the complexified K\"{a}hler form $\omega _{M}^{\mathbb{C}}=\omega
_{M}+i\beta _{M}$ is a complexified Calabi-Yau K\"{a}hler form if and only
if the complex valued function $\phi \left( x\right) +i\eta \left( x\right) $
satisfies the following complexified real Monge-Amper\'{e} equation, 
\begin{equation*}
\det \left( \phi _{jk}+i\eta _{jk}\right) =C,
\end{equation*}
for some nonzero constant $C$. If we write 
\begin{equation*}
\theta _{jk}\left( x\right) =\phi _{jk}\left( x\right) +i\eta _{jk}\left(
x\right) ,
\end{equation*}
then the above equation becomes $\det \left( \theta _{jk}\right) =C$. In
these notations, the complexified K\"{a}hler metric and complexified
K\"{a}hler form on $M$ are 
\begin{eqnarray*}
g_{M}^{\mathbb{C}} &=&\Sigma \theta _{jk}\left( x\right) \left(
dx^{j}\otimes dx^{k}+dy^{j}\otimes dy^{k}\right) \text{ and} \\
\omega _{M}^{\mathbb{C}} &=&\frac{i}{2}\Sigma \theta _{jk}\left( x\right)
dz^{j}\wedge d\bar{z}^{k},
\end{eqnarray*}
respectively.

Now we consider the dual $T^{n}$-invariant manifold $W$ as before. Instead
of the Legendre transformation $dx_{j}=\Sigma \phi _{jk}dx^{k}$, we need to
consider a complexified version of it. Symbolically we should write 
\begin{equation*}
dx_{j}=\Sigma \theta _{jk}dx^{k}=\Sigma \left( \phi _{jk}+i\eta _{jk}\right)
dx^{k}.
\end{equation*}
The precise meaning of this is the complex coordinates $dz_{j}$'s on $W$ is
determined by $\func{Re}dz_{j}=\phi _{jk}dx^{k}$ and $\func{Im}%
dz_{j}=dy_{j}+\eta _{jk}dx^{k}$. That is $dz_{j}=dx_{j}+idy_{j}$.

As before we define 
\begin{equation*}
\left( \theta ^{jk}\right) =\left( \theta _{jk}\right) ^{-1}.
\end{equation*}
It is easy to check directly that the canonical symplectic form on $W$ can
be expressed as follow 
\begin{eqnarray*}
\omega _{W}^{\mathbb{C}} &=&\Sigma dx^{j}\wedge dy_{j}, \\
&=&\frac{i}{2}\Sigma \theta ^{jk}dz_{j}\wedge d\bar{z}_{k}.
\end{eqnarray*}
Similarly the corresponding complexified K\"{a}hler metric is given by

\begin{eqnarray*}
g_{W}^{\mathbb{C}} &=&\Sigma \theta ^{jk}\left( dx_{j}\otimes
dx_{k}+dy_{j}\otimes dy_{k}\right) , \\
&=&\Sigma \theta ^{jk}dz_{j}\otimes d\bar{z}_{k}\text{.}
\end{eqnarray*}

After including the B-fields, we can argue using the same reasonings as
before and conclude: If we varies the complexified symplectic structures on $%
M$ while keeping its complex structure fixed. Then, under the Fourier
transformation along fibers and Legendre transformation on the base, it
corresponds to varying the complex structures on $W$ while keeping its
complexified symplectic structure fixed. And the reverse also hold true.

\begin{theorem}
Let $M$ be a $T^{n}$-invariant Calabi-Yau manifold and $W$ is its mirror.
Then the moduli space of complex structures on $M$ (resp. on $W$) is
identified with the moduli space of complexified symplectic structures on $W$
(resp. on $M$) under the above mirror transformation.
\end{theorem}

Remark: In order to have the above mirror transformation between complex
structures and symplectic structures, it is important that the B-fields
satisfy the complexified Monge-Amper\'{e} equation instead of being a
harmonic two form.

\bigskip

Second we use the second definition of a B-field, namely $\omega
^{n}=ci^{n}\Omega \bar{\Omega}$, $\func{Im}e^{i\theta }\left( \omega +i\beta
\right) ^{n}=0$ and $\func{Im}e^{i\phi }\Omega =0$ on the zero section. We
still use the Fourier and Legendre transformation, $\func{Re}dz_{j}=\phi
_{jk}dx^{k}$ and $\func{Im}dz_{j}=dy_{j}+\eta _{jk}dx^{k}$. Then Gross
observed that \cite{Gr}, 
\begin{eqnarray*}
\Omega _{W}\bar{\Omega}_{W} &=&\prod \left( \phi _{jk}dx^{k}+idy_{j}+i\eta
_{jk}dx^{k}\right) \left( \phi _{jk}dx^{k}-idy_{j}-i\eta _{jk}dx^{k}\right)
\\
&=&\prod \left( \phi _{jk}dx^{k}+idy_{j}\right) \left( \phi
_{jk}dx^{k}-idy_{j}\right) .
\end{eqnarray*}
So we still have $\omega _{W}^{n}=ci^{n}\Omega _{W}\bar{\Omega}_{W}$, as if $%
\beta $ has no effect.\ 

If we restrict $\Omega _{W}$ to the zero section of $W$, which is defined by 
$y_{j}=0$ for all $j$, then 
\begin{eqnarray*}
\func{Im}e^{i\theta }\Omega _{W} &=&\func{Im}e^{i\theta }\prod \left( \phi
_{jk}dx^{k}+idy_{j}+i\eta _{jk}dx^{k}\right) \\
&=&\func{Im}e^{i\theta }\prod \left( \phi _{jk}dx^{k}+i\eta
_{jk}dx^{k}\right) \\
&=&\func{Im}e^{i\theta }\det \left( \phi _{jk}+i\eta _{jk}\right)
dx^{1}\cdots dx^{n}\text{.}
\end{eqnarray*}
Hence the equation $\func{Im}e^{i\theta }\left( \omega +i\beta \right)
^{n}=0 $ for $\beta $ on the $M$ side is equivalent to the zero section of $%
W $ being a special Lagrangian submanifold $\func{Im}e^{i\theta }\Omega
_{W}=0$.

Hence under the mirror transformation, the following conditions on $M,$ 
\begin{eqnarray*}
\omega _{M}^{n} &=&ci^{n}\Omega _{M}\bar{\Omega}_{M} \\
\func{Im}e^{i\theta }\left( \omega _{M}+i\beta _{M}\right) ^{n} &=&0 \\
\func{Im}e^{i\phi }\Omega _{M} &=&0\text{ on the zero section.}
\end{eqnarray*}
becomes the corresponding conditions on $W$: 
\begin{eqnarray*}
\omega _{W}^{n} &=&ci^{n}\Omega _{W}\bar{\Omega}_{W} \\
\func{Im}e^{i\theta }\Omega _{W} &=&0\text{ on the zero section.} \\
\func{Im}e^{i\phi }\left( \omega _{W}+i\beta _{W}\right) ^{n} &=&0.
\end{eqnarray*}

In the following discussions, we will always use the first definition of the
B-field.

\bigskip

\section{Transforming $\Omega ^{p,q}$, $H^{p,q}\ $and Yukawa couplings}

\bigskip

\textbf{Transformation on moduli spaces: A holomorphic isometry}

Continue from above discussions, we are going to analyze the mirror
transformation from the moduli space of complexified symplectic structures
on $M$ to the moduli space of complex structures on $W$. We will see that
this map is both a holomorphic map and an isometry.

To do this, we need to know this transformation on the infinitesimal level.
Since infinitesimal deformation of K\"{a}hler structures on $M$ (resp.
complex structures on $W$) is parametrized by $H^{1}\left( M,T_{M}^{\ast
}\right) $ (resp. $H^{1}\left( W,T_{W}\right) $)\footnote{%
Cohomology groups are interpreted as spaces of $T^{n}$-invariant harmonic
forms.}, we should have a homomorphism 
\begin{equation*}
T:H^{1}\left( M,T_{M}^{\ast }\right) \rightarrow H^{1}\left( W,T_{W}\right) 
\text{.}
\end{equation*}
Suppose we vary the $T^{n}$ symplectic form on $M$ to 
\begin{equation*}
\omega _{M}^{new}=\omega _{M}+\varepsilon \Sigma \xi
_{jk}dx^{j}dy^{k}=\Sigma \left( \phi _{jk}+\varepsilon \xi _{jk}\right)
dx^{j}dy^{k}.
\end{equation*}
Here 
\begin{equation*}
\xi _{jk}=\frac{\partial ^{2}\xi \left( x\right) }{\partial x^{j}\partial
x^{k}}\text{,}
\end{equation*}
and $\varepsilon $ is the deformation parameter. Then 
\begin{equation*}
\Sigma \xi _{jk}dx^{j}dy^{k}=\frac{i}{2}\Sigma \xi _{jk}dz^{j}d\bar{z}^{k}
\end{equation*}
represents an element in $\Omega ^{0,1}\left( M,T_{M}^{\ast }\right) $ which
parametrizes deformations of K\"{a}hler forms. This form is harmonic, namely
it defines an element in $H^{1}\left( M,T_{M}^{\ast }\right) $, if and only
if $\Sigma _{k}\xi _{jkk}=0$ for all $j$. If we assume every member of the
family of $T^{n}$-invariant K\"{a}hler forms is Calabi-Yau, then the
infinitesimal variation $\xi $ satisfies a linearization of the Monge-Amper%
\'{e} equation. This implies that $\xi $ is harmonic.

Then the new complex structure on $W$ is determined by its new complex
coordinates 
\begin{eqnarray*}
dz_{j}^{new} &=&\Sigma \left( \phi _{jk}+\varepsilon \xi _{jk}\right)
dx^{k}+idy_{j} \\
&=&dz_{j}+\varepsilon \Sigma \xi _{jk}dx^{k} \\
&=&\Sigma \left( \delta _{j}^{l}+\frac{\varepsilon }{2}\phi ^{lk}\xi
_{jk}\right) dz_{l}+\frac{\varepsilon }{2}\phi ^{kl}\xi _{jk}d\bar{z}_{l}.
\end{eqnarray*}
Therefore if we project the new $\bar{\partial}$-operator on $W$ to the old $%
\Omega ^{0,1}\left( W\right) $, we have 
\begin{equation*}
\bar{\partial}^{new}=\bar{\partial}-\frac{\varepsilon }{2}\Sigma \phi
^{jk}\xi _{kl}\frac{\partial }{\partial z_{l}}\otimes d\bar{z}_{j}+O\left(
\varepsilon ^{2}\right) \text{.}
\end{equation*}
It gives an element 
\begin{equation*}
-\frac{1}{2}\Sigma \phi ^{jk}\xi _{kl}\frac{\partial }{\partial z_{l}}%
\otimes d\bar{z}_{j}\in \Omega ^{0,1}\left( W,T_{W}\right) ,
\end{equation*}
that determines the infinitesimal deformation of corresponding complex
structures on $W$. This element is harmonic, namely it defines an element in 
$H^{1}\left( W,T_{W}\right) $, if and only if $\frac{\partial }{\partial
x_{j}}\left( \xi _{lk}\phi ^{jl}\right) =0$. This is equivalent to $\xi
_{jjk}=0$. Hence we have obtained explicitly the homomorphism 
\begin{eqnarray*}
H^{1}\left( M,T_{M}^{\ast }\right)  &\rightarrow &H^{1}\left( W,T_{W}\right) 
\\
i\Sigma \xi _{jk}dz^{j}d\bar{z}^{k} &\rightarrow &-\Sigma \xi _{jk}\phi ^{kl}%
\frac{\partial }{\partial z_{j}}\otimes d\bar{z}_{l}.
\end{eqnarray*}
Notice that these infinitesimal deformations are $T^{n}$-invariant, $\xi
_{jk}=\xi _{jk}\left( x\right) $. Therefore we can use $\xi _{jk}$'s to
denote both a tensor in $M$ and its transformation in $W$.

We should also include the B-fields and use the complexified symplectic
forms on $M$, however the formula is going to be the same (with $\theta $
replacing $\phi $). From this description, it is obvious that the
transformation from the moduli space of complexified symplectic forms on $M$
to the moduli space of complex structures on $W$ is holomorphic.

Next we are going to verify that this mirror map between the two moduli
spaces is an isometry. We take two such deformation directions $i\Sigma \xi
_{jk}dz^{j}d\bar{z}^{k}$ and $i\Sigma \zeta _{jk}dz^{j}d\bar{z}^{k}$, their L%
$^{2}$-inner product is given by 
\begin{equation*}
\left\langle i\Sigma \xi _{jk}dz^{j}d\bar{z}^{k},i\Sigma \zeta _{jk}dz^{j}d%
\bar{z}^{k}\right\rangle _{M}=2V\int_{D}\phi ^{jl}\phi ^{km}\xi _{jk}\zeta
_{lk}dv_{D}.
\end{equation*}
While the L$^{2}$-inner product of their image on the $W$ side is given by 
\begin{eqnarray*}
&&\left\langle -\Sigma \xi _{jk}\phi ^{kl}\frac{\partial }{\partial z_{j}}%
\otimes d\bar{z}_{l},-\Sigma \zeta _{jk}\phi ^{kl}\frac{\partial }{\partial
z_{j}}\otimes d\bar{z}_{l}\right\rangle _{W} \\
&=&2V^{-1}\int_{D}\phi ^{jp}\phi _{lq}\left( \zeta _{jk}\phi ^{kl}\xi
_{pm}\phi ^{mq}\right) dv_{D} \\
&=&2V^{-1}\int_{D}\phi ^{jl}\phi ^{km}\xi _{jk}\zeta _{lk}dv_{D}.
\end{eqnarray*}
Here $V$ (resp. $V^{-1}$) is the volume of the special Lagrangian fiber in $%
M $ (resp. $W$). Therefore up to a overall constant, this transformation
between the two moduli spaces is not just holomorphic, it is an isometry
too. We conclude that

\begin{theorem}
The above explicit mirror map from the moduli space of complex structures on 
$M$ (resp. on $W$) to the moduli space of complexified symplectic structures
on $W$ (resp. on $M$) is a holomorphic isometry.
\end{theorem}

\bigskip

\textbf{Transforming differential forms}

Next we transform differential forms of higher degrees from $M$ to $W$: 
\begin{equation*}
T:\Omega ^{0,q}\left( M,\Lambda ^{p}T_{M}^{\ast }\right) \rightarrow \Omega
^{0,q}\left( W,\Lambda ^{p}T_{W}\right) .
\end{equation*}
Using the triviality of the canonical line bundle of $W$, this is the same
as 
\begin{equation*}
T:\Omega ^{p,q}\left( M\right) \rightarrow \Omega ^{n-p,q}\left( W\right) 
\text{.}
\end{equation*}
Readers are reminded that we are discussing only $T^{n}$-invariant
differential forms.

\bigskip

First we give the motivations for this homomorphism. Since $M$ and $W$ are
related by fiberwise dual torus construction, the obvious transformation for
their tensors would be 
\begin{equation*}
\left\{ 
\begin{array}{l}
\medskip dx^{j}\rightarrow dx^{j}=\Sigma \phi ^{jk}dx_{k} \\ 
dy^{j}\rightarrow \frac{\partial }{\partial y_{j}}\text{.}
\end{array}
\right.
\end{equation*}
In symplectic language, such transformation uses the real polarizations of $%
M $ and $W$. To transform $\left( p,q\right) $ forms, we want to map this
real polarization to the complex polarization. The real polarization is
defined by the vertical tangent bundle $V\subset T_{M}$ and the complex
polarization is defined by $T_{M}^{1,0}\subset T_{M}\otimes \mathbb{C}$. So
it is natural to carry $V\otimes \mathbb{C}$ to $T_{M}^{1,0}$ and its
complement to $T_{M}^{0,1}$. That is $dz^{j}\rightarrow dy^{j}$ and $d\bar{z}%
^{j}\rightarrow dx^{j}$ on the $M$ side. By doing the same identification on
the $W$ side and compose with the above transformation, we have 
\begin{equation*}
T:\Omega ^{0,q}\left( M,\Lambda ^{p}T_{M}^{\ast }\right) \rightarrow \Omega
^{0,q}\left( W,\Lambda ^{p}T_{W}\right) ,
\end{equation*}
with 
\begin{eqnarray*}
T\left( dz^{j}\right) &=&\frac{\partial }{\partial z_{j}} \\
T\left( d\bar{z}^{j}\right) &=&\Sigma \phi ^{jk}d\bar{z}_{k}\text{.}
\end{eqnarray*}

This homomorphism obviously coincides with the previous identification
between infinitesimal deformation of symplectic structures on $M$ and
complex structures on $W$ up to a constant factor $i$.

Using the holomorphic volume form $\Omega _{W}=dz_{1}dz_{2}\cdots dz_{n}$ on 
$W$, we can identify $\wedge ^{p}T_{W}^{1,0}$ with $\Lambda
^{n-p}T_{W}^{1,0\ast }$, so we obtain a homomorphism 
\begin{equation*}
T:\Omega ^{p,q}\left( M\right) \rightarrow \Omega ^{n-p,q}\left( W\right) 
\text{.}
\end{equation*}
Explicitly if 
\begin{equation*}
\alpha =\Sigma \alpha _{i_{1}...i_{p}\bar{j}_{1}...\bar{j}%
_{q}}dz^{i_{1}}\cdots dz^{i_{p}}d\bar{z}^{j_{1}}\cdots d\bar{z}^{j_{q}}\in
\Omega ^{p,q}\left( M\right)
\end{equation*}
then 
\begin{equation*}
T\left( \alpha \right) =\Sigma \alpha _{i_{1}...i_{p}\bar{j}_{1}...\bar{j}%
_{q}}\phi ^{k_{1}\bar{j}_{1}}\cdots \phi ^{k_{q}\bar{j}_{q}}dz_{1}\cdots 
\widehat{dz_{i_{1}}}\cdots \widehat{dz_{i_{p}}}\cdots dz_{n}d\bar{z}%
_{k_{1}}\cdots d\bar{z}_{k_{q}}.
\end{equation*}

\bigskip

\bigskip

\textbf{Transforming }$\mathbf{H}^{p,q}\left( M\right) $\textbf{\ to }$%
\mathbf{H}^{n-p,q}\left( W\right) $

If $\alpha \in \Omega ^{p,q}\left( M\right) $ is a $T^{n}$-invariant form,
then we claim that the above transformation of differential forms from $M$
to $W$ does commute with the $\bar{\partial}$-operator and also $\bar{%
\partial}^{\ast }$-operator. Therefore it descends to the Hodge cohomology
(also Dolbeault cohomology) level: 
\begin{equation*}
T:H^{p,q}\left( M\right) \rightarrow H^{n-p,q}\left( W\right) \text{.}
\end{equation*}
To simplify our notations we assume that $\alpha $ is of type $\left(
1,1\right) .$ That is 
\begin{equation*}
\alpha =\Sigma \alpha _{jk}dz^{j}\wedge d\bar{z}^{k}.
\end{equation*}
The form $\alpha $ being $T^{n}$-invariant means that $\alpha _{jk}=\alpha
_{jk}\left( x\right) $ depends on the $x$ variables only. We have 
\begin{equation*}
\bar{\partial}\alpha =\frac{1}{2}\Sigma \left( \frac{\partial \alpha _{jk}}{%
\partial x^{p}}-\frac{\partial \alpha _{jp}}{\partial x^{k}}\right) dz^{j}d%
\bar{z}^{p}d\bar{z}^{k}
\end{equation*}
Their transformations are 
\begin{eqnarray*}
T\left( \alpha \right)  &=&\Sigma \alpha _{jk}dz_{1}\cdots \widehat{dz_{j}}%
\cdots dz_{n}d\bar{z}_{l}\phi ^{kl}, \\
T\left( \bar{\partial}\alpha \right)  &=&\frac{1}{2}\Sigma \left( \frac{%
\partial \alpha _{jk}}{\partial x^{p}}-\frac{\partial \alpha _{jp}}{\partial
x^{k}}\right) dz_{1}\cdots \widehat{dz_{j}}\cdots dz_{n}d\bar{z}_{q}d\bar{z}%
_{l}\phi ^{kl}\phi ^{pq}.
\end{eqnarray*}

Now 
\begin{equation*}
\bar{\partial}T\left( \alpha \right) =\frac{1}{2}\Sigma \left( \frac{%
\partial }{\partial x_{q}}\left( \alpha _{jk}\phi ^{kl}\right) -\frac{%
\partial }{\partial x_{l}}\left( \alpha _{jk}\phi ^{kq}\right) \right)
dz_{1}\cdots \widehat{dz_{j}}\cdots dz_{n}d\bar{z}_{q}d\bar{z}_{l}.
\end{equation*}
Using the Legendre transformation 
\begin{equation*}
\frac{\partial }{\partial x_{q}}=\Sigma \phi ^{pq}\frac{\partial }{\partial
x^{p}},
\end{equation*}
we have 
\begin{eqnarray*}
&&\frac{\partial }{\partial x_{q}}\left( \alpha _{jk}\phi ^{kl}\right) -%
\frac{\partial }{\partial x_{l}}\left( \alpha _{jk}\phi ^{kq}\right) \\
&=&\Sigma \phi ^{pq}\frac{\partial \alpha _{jk}}{\partial x^{p}}\phi
^{kl}+\phi ^{pq}\alpha _{jk}\frac{\partial \phi ^{kl}}{\partial x^{p}}-\phi
^{pl}\frac{\partial \alpha _{jk}}{\partial x^{p}}\phi ^{kq}-\phi ^{pl}\alpha
_{jk}\frac{\partial \phi ^{kq}}{\partial x^{p}} \\
&=&\Sigma \left( \frac{\partial \alpha _{jk}}{\partial x^{p}}-\frac{\partial
\alpha _{jp}}{\partial x^{k}}\right) \phi ^{kl}\phi ^{pq}+\alpha _{jk}\left(
-\phi ^{pq}\phi ^{ks}\frac{\partial \phi _{st}}{\partial x^{p}}\phi
^{tl}+\phi ^{pl}\phi ^{ks}\frac{\partial \phi _{st}}{\partial x^{p}}\phi
^{tq}\right) .
\end{eqnarray*}
The second bracket vanishes because $\frac{\partial }{\partial x^{p}}\phi
_{st}$ is symmetric with respect to $s,t$ and $p$. Hence we have 
\begin{equation*}
\bar{\partial}T\left( \alpha \right) =T\left( \bar{\partial}\alpha \right) 
\text{.}
\end{equation*}
We can also verify 
\begin{equation*}
\bar{\partial}^{\ast }T\left( \alpha \right) =T\left( \bar{\partial}^{\ast
}\alpha \right)
\end{equation*}
in the same way, and is left to our readers. Therefore the transformation $T$
descends to both the Hodge cohomology and Dolbeault cohomology. Moreover if
we go the other direction, namely from $W$ to $M$, then the corresponding
transformation is the inverse of $T$. That is $T$ is an isomorphism. So we
have the following result.

\begin{theorem}
The above mirror transformation $T$ identifies $\Omega ^{p,q}\left( M\right) 
$ (resp. $H^{p,q}\left( M\right) $) with $\Omega ^{n-p,q}\left( W\right) $
(resp. $H^{n-p,q}\left( W\right) $).
\end{theorem}

\bigskip

\textbf{Transforming Yukawa couplings}

Next we compare Yukawa couplings on these moduli spaces of complex and
symplectic structures; they are first computed by Mark Gross in \cite{Gr}.
We choose any $n$ closed differential forms of type $\left( 1,1\right) $ on $%
M$: $\alpha ,\beta ,\cdots ,\gamma $. We write $\alpha =\Sigma \alpha
_{ij}dz^{i}\wedge d\bar{z}^{j}$ and so on. The Yukawa coupling in the A side
on $M$ is defined and computed as follows, 
\begin{eqnarray*}
_{A}Y_{M}\left( \alpha ,\beta ,...,\gamma \right)  &=&\int_{M}\alpha \wedge
\beta \wedge \cdots \wedge \gamma  \\
&=&\int_{M}\pm \alpha _{i_{1}j_{1}}\beta _{i_{2}j_{2}}\cdots \gamma
_{i_{n}j_{n}}dv_{M} \\
&=&V\int_{x}\sum \pm \alpha _{i_{1}j_{1}}\beta _{i_{2}j_{2}}\cdots \gamma
_{i_{n}j_{n}}dx^{1}dx^{2}\cdots dx^{n}
\end{eqnarray*}
where the summation is such that $\left\{ i_{1},i_{2},...,i_{n}\right\}
=\left\{ j_{1},j_{2},...,j_{n}\right\} =\left\{ 1,2,...,n\right\} .$ The
constant $V$ is the volume of a special Lagrangian fiber in $M$.

For the Yukawa coupling in the $B$ side on $W$ we have the following
definition: For $\alpha ^{\prime }=T\left( \alpha \right) ,\beta ^{\prime
}=T\left( \beta \right) ,\cdots ,\gamma ^{\prime }=T\left( \gamma \right)
\in \Omega ^{0,1}\left( W,T_{W}\right) $, we have

\begin{equation*}
_{B}Y_{W}\left( \alpha ^{\prime },\beta ^{\prime },...,\gamma ^{\prime
}\right) =\int_{W}\Omega \wedge \delta _{\alpha ^{\prime }}\delta _{\beta
^{\prime }}...\delta _{\gamma ^{\prime }}\Omega .
\end{equation*}
Since $\Omega =dz_{1}dz_{2}\cdots dz_{n}$ with $dz_{j}=\Sigma \phi
_{jk}dx^{k}+idy_{j}$, we have 
\begin{equation*}
\delta _{\alpha }\Omega =\Sigma \alpha
_{1k}dx^{k}dz_{2}...dz_{n}+dz_{1}\alpha
_{2k}dx^{k}...dz_{n}+...+dz_{1}dz_{2}...\alpha _{nk}dx^{k}.
\end{equation*}
Similarly we obtain 
\begin{equation*}
\delta _{\alpha ^{\prime }}\delta _{\beta ^{\prime }}...\delta _{\gamma
^{\prime }}\Omega =\sum \pm \alpha _{i_{1}j_{1}}\beta _{i_{2}j_{2}}\cdots
\gamma _{i_{n}j_{n}}dx^{1}dx^{2}\cdots dx^{n}
\end{equation*}
and therefore, up to an overall constant, we have 
\begin{equation*}
_{A}Y_{M}\left( \alpha ,\beta ,...,\gamma \right) =_{B}Y_{W}\left( \alpha
^{\prime },\beta ^{\prime },...,\gamma ^{\prime }\right) .
\end{equation*}

\begin{theorem}
\cite{Gr} The above mirror transformation identifies the Yukawa coupling on
the moduli spaces of complexified symplectic structures on $M$ (resp. on $W$%
) with the Yukawa coupling on the moduli space of complex structures on $W$
(resp. on $M$).
\end{theorem}

Remark: The Yukawa coupling is the $n^{th}$ derivative of a local
holomorphic function on the moduli space, called the prepotential $\mathcal{F%
}$. On the A-side, this is given by 
\begin{equation*}
_{A}\mathcal{F}\left( M\right) =\int_{M}\omega ^{n}.
\end{equation*}
On the B-side, we need to specify a holomorphic family of the holomorphic
volume form locally on the moduli space of complex structure on $W$. Then
the prepotential function is given by 
\begin{equation*}
_{B}\mathcal{F}\left( W\right) =\int_{W}\Omega \wedge \bar{\Omega}.
\end{equation*}
Similarly we can identify these two prepotentials by this transformation.

In fact we can express this identification of the two moduli spaces,
together with identifications of all these structures on them, namely $%
\Omega ^{\ast ,\ast },H^{\ast ,\ast },\mathcal{Y}$ and $\mathcal{F}$, as an
isomorphism of two Frobenius manifolds.

\section{$\mathbf{sl}_{2}\mathbf{\times sl}_{2}$-action on cohomology and
their mirror transform}

In this section we show that on the levels of differential forms and
cohomology of $M$, there are two commuting $\mathbf{sl}\left( 2\right) $ Lie
algebra actions. Moreover the mirror transformation between $M$ and $W$
interchanges them. The first $\mathbf{sl}\left( 2\right) $ action exists for
all K\"{a}hler manifolds. We should note that the results of this section
depends only on the $T^{n}$-invariant condition but not the Calabi-Yau
condition.

This type of structure is first proposed by Gopakumar and Vafa in \cite{GV2}
on the moduli space of flat $U\left( 1\right) $ bundles over curves in $M$,
that is B-cycles. They conjectured that this $\mathbf{sl}\left( 2\right)
\times \mathbf{sl}\left( 2\right) $ representation determines all
Gromov-Witten invariants in every genus in a Calabi-Yau manifold. In fact we
conjecture that such $\mathbf{sl}\left( 2\right) \times \mathbf{sl}\left(
2\right) $ action on cohomology groups should exist for every moduli space
of A- or B-cycles (with the rank of the bundle equals one) on mirror
manifolds $M$ and $W$ \cite{L1}.

\bigskip

\textbf{Hard Lefschetz }$\mathbf{sl}\left( 2\right) $ \textbf{action}

Recall that the cohomology of any K\"{a}hler manifold admits an $\mathbf{sl}%
\left( 2\right) $ action. Let us recall its construction: Let $M$ be a
K\"{a}hler manifold with K\"{a}hler form $\omega _{M}$. Wedging with $\omega
_{M}$ gives a homomorphism $L_{A}$: 
\begin{equation*}
L_{A}:\Omega ^{k}\left( M\right) \rightarrow \Omega ^{k+2}\left( M\right) .
\end{equation*}
Let 
\begin{equation*}
\Lambda _{A}:\Omega ^{k+2}\left( M\right) \rightarrow \Omega ^{k}\left(
M\right)
\end{equation*}
be its adjoint homomorphism. Then we have the following relations 
\begin{equation*}
\left[ L_{A},\Lambda _{A}\right] =H_{A},
\end{equation*}
where $H_{A}=\left( n-k\right) I$ is the multiplication endomorphism on $%
\Omega ^{k}\left( M\right) $. Moreover we have 
\begin{eqnarray*}
\left[ L_{A},H_{A}\right] &=&\,2L_{A}, \\
\left[ \Lambda _{A},H_{A}\right] &=&-2\Lambda _{A}.
\end{eqnarray*}
These commutating relations determine an $\mathbf{sl}\left( 2\right) $
action on $\Omega ^{\ast }\left( M\right) $. We call it the hard Lefschetz $%
\mathbf{sl}\left( 2\right) $ action.

These operations commute with $\bar{\partial}$ and $\bar{\partial}^{\ast }$
because $\omega _{M}$ is a parallel form on $M$. Therefore this $\mathbf{sl}%
\left( 2\right) $ action descends to the cohomology group $H^{\ast ,\ast
}\left( M\right) $.

\bigskip

\textbf{Variation of Hodge structures }$sl\left( 2\right) $\textbf{\ action }

Suppose $M$ is a $T^{n}$-invariant manifold. It comes with a natural family
of deformation of complex structures whose complex coordinates are given by $%
z^{j}=\frac{1}{t}x^{j}+iy^{j}$.

Recall from the standard deformation theory that a deformation of complex
structures determines a variation of Hodge structures. Infinitesimally the
variation of Hodge filtration $F^{p}\left( H^{\ast }\left( M,\mathbb{C}%
\right) \right) $ lies insides $F^{p-1}\left( H^{\ast }\left( M,\mathbb{C}%
\right) \right) $: If we write the infinitesimal variation of complex
structure as $\frac{dM_{t}}{dt}\in H^{1}\left( M,T_{M}\right) $. Then the
variation of Hodge structures is determined by taking the trace of the cup
product with $\frac{dM_{t}}{dt}$ which sends $H^{q}\left( M,\Omega
_{M}^{p}\right) $ to $H^{q+1}\left( M,\Omega _{M}^{p-1}\right) $. We denote
this homomorphism by $L_{B}$. That is 
\begin{equation*}
L_{B}=\frac{dM_{t}}{dt}:H^{p,q}\left( M\right) \rightarrow H^{p-1,q+1}\left(
M\right) \text{.}
\end{equation*}
For the $T^{n}$-invariant K\"{a}hler manifold $M$, it turns out that $L_{B}$
determines an $\mathbf{sl}\left( 2\right) $ action on $H^{\ast }\left( M,%
\mathbb{C}\right) =\oplus H^{p,q}\left( M\right) $.

To describe this $\mathbf{sl}\left( 2\right) $ action explicitly, first we
need to describe the adjoint of $L_{B}$ which we will call $\Lambda _{B}$.
In general if 
\begin{equation*}
\frac{dM_{t}}{dt}=\Sigma a_{\bar{k}}^{j}\left( z,\bar{z}\right) \frac{%
\partial }{\partial z^{j}}\otimes d\bar{z}^{k}
\end{equation*}
on a K\"{a}hler manifold with metric $\Sigma g_{j\bar{k}}dz^{j}\otimes d\bar{%
z}^{k}$, then the adjoint of $L_{B}$ on the level of differential forms is
just 
\begin{equation*}
\Lambda _{B}=\Sigma b_{j}^{\bar{k}}\frac{\partial }{\partial \bar{z}^{k}}%
\otimes dz^{j}
\end{equation*}
where $b_{j}^{\bar{k}}=\overline{g^{k\bar{l}}a_{\bar{l}}^{m}g_{m\bar{j}}}$.
Since $\bar{\partial}\left( \frac{dM_{t}}{dt}\right) =0$, $L_{B}$ commutes
with the $\bar{\partial}$-operator. However $\Lambda _{B}$ might not
commutes with $\bar{\partial}$ and therefore would not descend to the level
of cohomology in general.

For the $T^{n}$-invariant case, it is not difficult to check directly that $%
\frac{dM_{t}}{dt}=\frac{t-1}{2}\Sigma \frac{\partial }{\partial z^{j}}%
\otimes d\bar{z}^{j}$. That is the whole family of complex structures on $M$
is along the same direction. We can rescale and assume 
\begin{equation*}
\frac{dM_{t}}{dt}=\Sigma \frac{\partial }{\partial z^{j}}\otimes d\bar{z}%
^{j}.
\end{equation*}
That is $a_{\bar{k}}^{j}\left( z,\bar{z}\right) =\delta _{jk}$. Hence 
\begin{eqnarray*}
b_{j}^{\bar{k}} &=&\overline{g^{k\bar{l}}a_{\bar{l}}^{m}g_{m\bar{j}}} \\
&=&\phi ^{kl}\delta _{ml}\phi _{mj} \\
&=&\delta _{jk}\text{.}
\end{eqnarray*}
That is 
\begin{equation*}
\Lambda _{B}=\Sigma \frac{\partial }{\partial \bar{z}^{j}}\otimes dz^{j}.
\end{equation*}

Moreover their commutator $H_{B}=\left[ L_{B},\Lambda _{B}\right] $ is the
multiplication of $\left( p-q\right) $ on forms in $\Omega ^{p,q}\left(
M\right) $. We have the following result:

\begin{theorem}
On a $T^{n}$-invariant manifold $M$ as before, if we define 
\begin{eqnarray*}
L_{B} &=&\Sigma \frac{\partial }{\partial z^{j}}\otimes d\bar{z}^{j}:\Omega
^{p,q}\left( M\right) \rightarrow \Omega ^{p-1,q+1}\left( M\right)  \\
\Lambda _{B} &=&\Sigma \frac{\partial }{\partial \bar{z}^{j}}\otimes
dz^{j}:\Omega ^{p,q}\left( M\right) \rightarrow \Omega ^{p+1,q-1}\left(
M\right)  \\
H_{B} &=&\left[ L_{B},\Lambda _{B}\right] =\left( p-q\right) :\Omega
^{p,q}\left( M\right) \rightarrow \Omega ^{p,q}\left( M\right) .
\end{eqnarray*}
Then they satisfy $\Lambda _{B}=\left( L_{B}\right) ^{\ast }$ and 
\begin{equation*}
\begin{array}{ccc}
\left[ L_{B},\Lambda _{B}\right]  & = & H_{B} \\ 
\left[ H_{B},L_{B}\right]  & = & -2L_{B} \\ 
\left[ H_{B},\Lambda _{B}\right]  & = & \,2\Lambda _{B}\text{.}
\end{array}
\end{equation*}
Hence they define an $\mathbf{sl}\left( 2\right) $ action on $\Omega ^{\ast
,\ast }\left( M\right) $. Moreover these operators commute with $\bar{%
\partial}$ and $\bar{\partial}^{\ast }$ and descend to give an $\mathbf{sl}%
\left( 2\right) $ action on $H^{\ast }\left( M,\mathbb{C}\right) $.
\end{theorem}

As a corollary we have the following.

\begin{corollary}
On any $T^{n}$-invariant manifold $M$ with $\frac{dM_{t}}{dt}$ as before.
The operators $L_{B}$ defined by the variation of Hodge structures, its
adjoint operator $\Lambda _{B}$ and their commutator $H_{B}=\left[
L_{B},\Lambda _{B}\right] $ together defines an $\mathbf{sl}\left( 2\right) $
action on the cohomology of $M$.
\end{corollary}

We call this the variation of Hodge structures $\mathbf{sl}\left( 2\right) $
action, or simply VHS $\mathbf{sl}\left( 2\right) $ action.

\bigskip

\textbf{An }$\mathbf{sl}\left( 2\right) \times \mathbf{sl}\left( 2\right) $ 
\textbf{action on cohomology}

We already have two $\mathbf{sl}\left( 2\right) $ actions on $H^{\ast
}\left( M\right) $, we want to show that they commute with each other.

\begin{lemma}
On a $T^{n}$-invariant manifold $M$ as above, we have 
\begin{eqnarray*}
\left[ L_{A},L_{B}\right]  &=&0, \\
\left[ L_{A},\Lambda _{B}\right]  &=&0.
\end{eqnarray*}
\end{lemma}

Proof of lemma: We verify this lemma by direct calculations. Let us consider 
\begin{eqnarray*}
&&L_{A}L_{B}\left( dz^{j_{1}}\cdots dz^{j_{p}}d\bar{z}^{k_{1}}\cdots d\bar{z}%
^{k_{q}}\right) \\
&=&L_{A}\Sigma \left( -1\right) ^{p-s}dz^{j_{1}}\cdots \widehat{dz^{j_{s}}}%
\cdots dz^{j_{p}}d\bar{z}^{j_{s}}d\bar{z}^{k_{1}}\cdots d\bar{z}^{k_{q}} \\
&=&\Sigma \left( -1\right) ^{p-s}\left( -1\right) ^{p-1}\phi
_{jk}dz^{j}dz^{j_{1}}\cdots \widehat{dz^{j_{s}}}\cdots dz^{j_{p}}d\bar{z}%
^{k}d\bar{z}^{j_{p}}d\bar{z}^{k_{1}}\cdots d\bar{z}^{k_{q}}\text{.}
\end{eqnarray*}
On the other hand, 
\begin{eqnarray*}
&&L_{B}L_{A}\left( dz^{j_{1}}\cdots dz^{j_{p}}d\bar{z}^{k_{1}}\cdots d\bar{z}%
^{k_{q}}\right) \\
&=&L_{B}\Sigma \left( -1\right) ^{p}\phi _{jk}dz^{j}dz^{j_{1}}\cdots
dz^{j_{p}}d\bar{z}^{k}d\bar{z}^{k_{1}}\cdots d\bar{z}^{k_{q}}.
\end{eqnarray*}
If $j$ is not any of the $j_{r}$'s, then 
\begin{eqnarray*}
&&L_{B}\left( -1\right) ^{p}\phi _{jk}dz^{j}dz^{j_{1}}\cdots dz^{j_{p}}d\bar{%
z}^{k}d\bar{z}^{k_{1}}\cdots d\bar{z}^{k_{q}} \\
&=&\Sigma \phi _{jk}dz^{j_{1}}\cdots dz^{j_{p}}d\bar{z}^{j}d\bar{z}%
^{k}dz^{k_{1}}\cdots dz^{k_{q}} \\
&&+\Sigma \left( -1\right) ^{p}\left( -1\right) ^{p-s}\phi
_{jk}dz^{j}dz^{j_{1}}\cdots \widehat{dz^{j_{s}}}\cdots dz^{j_{p}}d\bar{z}%
^{j_{s}}d\bar{z}^{k}d\bar{z}^{k_{1}}\cdots d\bar{z}^{k_{q}}.
\end{eqnarray*}
However the first term on the right hand side is zero because $\phi
_{jk}=\phi _{kj}$ and $d\bar{z}^{j}d\bar{z}^{k}=-d\bar{z}^{k}d\bar{z}^{j}$.
If $j$ is one of the $j_{r}$'s, it turns out we have the same result. This
verifies $L_{A}L_{B}=L_{B}L_{A}$ on such forms. However forms of this type
generate all differential form and therefore we have 
\begin{equation*}
\left[ L_{A},L_{B}\right] =0.
\end{equation*}

If we replace $L_{B}=\Sigma \frac{\partial }{\partial z^{j}}\otimes d\bar{z}%
^{j}$ by $\Lambda _{B}=\Sigma \frac{\partial }{\partial \bar{z}^{j}}\otimes
dz^{j}$, it is not difficult to check that the same argument works and give
us 
\begin{equation*}
\left[ L_{A},\Lambda _{B}\right] =0\text{.}
\end{equation*}
Hence we have the lemma. $\square $

\begin{corollary}
On the cohomology of $M$ as above, the hard Lefschetz $\mathbf{sl}\left(
2\right) $ action and the VHS $\mathbf{sl}\left( 2\right) $ action commute.

In other words, we have an $\mathbf{sl}\left( 2\right) \mathbf{\times sl}%
\left( 2\right) $ action on $H^{\ast }\left( M,\mathbb{C}\right) $.
\end{corollary}

Proof of corollary: From the lemma we have $\left[ L_{A},L_{B}\right] =0$
and $\left[ L_{A},\Lambda _{B}\right] =0.$ Taking adjoint, we obtain the
other commutation relations. Hence the result. $\square $

\bigskip

Remark: The hard Lefschetz $\mathbf{sl}\left( 2\right) $ action is a
vertical action and the one from the variation of Hodge structure is a
horizontal action with respect to the Hodge diamond in the following sense: $%
L_{A}\left( H^{p,q}\right) \subset H^{p+1,q+1}$ and $L_{B}\left(
H^{p,q}\right) \subset H^{p-1,q+1}$.

Remark: Notice that $\mathbf{so}\left( 3,1\right) \mathbf{=sl}\left(
2\right) \mathbf{\times sl}\left( 2\right) $. Later we will show that when $%
M $ is a hyperk\"{a}hler manifold, this $\mathbf{so}\left( 3,1\right) $
action embeds naturally inside the canonical hyperk\"{a}hler $\mathbf{so}%
\left( 4,1\right) $ action on its cohomology group.

\bigskip

\textbf{Transforming the }$\mathbf{sl}\left( 2\right) \mathbf{\times sl}%
\left( 2\right) $\textbf{\ action}

First we recall that the variation of complex structures $dz^{j}=\frac{1}{t}%
dx^{j}+idy^{j}$ on $M$\ was carried to the variation of symplectic
structures $\omega =\frac{1}{t}dx^{j}dy_{j}$\ on $W$.

\begin{theorem}
Let $M$ and $W$ be mirror $T^{n}$-invariant K\"{a}hler manifolds to each
other. Then the mirror transformation $T$ carries the hard Lefschetz $%
\mathbf{sl}\left( 2\right) $ action on $M$ (resp. on $W$) to the variation
of Hodge structure $\mathbf{sl}\left( 2\right) $ action on $W$ (resp. on $M$%
).
\end{theorem}

Proof: Let us start by comparing $H_{A}$ and $H_{B}$. On $\Omega
^{p,q}\left( M\right) $, $H_{B}$ is the multiplication by $p-q$. On $\Omega
^{n-p,q}\left( W\right) $, $H_{A}$ is the multiplication by $n-\left(
n-p\right) -q=p-q$. On the other hand, $T$ carries $\Omega ^{p,q}\left(
M\right) $ to $\Omega ^{n-p,q}\left( W\right) $. Therefore 
\begin{equation*}
H_{A}T=TH_{B}\text{.}
\end{equation*}

Next we compare $L_{A}$ and $L_{B}$ for one forms on $M$. For the $\left(
0,1\right) $ form $d\bar{z}^{j}$, we have 
\begin{eqnarray*}
T\left( d\bar{z}^{j}\right) &=&\Sigma \phi ^{jk}dz_{1}\cdots dz_{n}d\bar{z}%
_{k}, \\
L_{A}T\left( d\bar{z}^{j}\right) &=&0.
\end{eqnarray*}
The last equality follows from type considerations. On the other hand $%
L_{B}\left( d\bar{z}^{j}\right) =0$, therefore 
\begin{equation*}
L_{A}T\left( d\bar{z}^{j}\right) =TL_{B}\left( d\bar{z}^{j}\right) \text{.}
\end{equation*}

For the $\left( 1,0\right) $ form $dz^{j}$, we have 
\begin{eqnarray*}
T\left( dz^{j}\right) &=&\left( -1\right) ^{n-j}dz_{1}\cdots \widehat{dz_{j}}%
\cdots dz_{n} \\
L_{A}T\left( dz^{j}\right) &=&\Sigma \left( -1\right) ^{n-j}\left( -1\right)
^{n+j}\phi ^{jk}dz_{1}\cdots dz_{n}d\bar{z}_{k} \\
&=&\Sigma \phi ^{jk}dz_{1}\cdots dz_{n}d\bar{z}_{k}.
\end{eqnarray*}
On the other hand 
\begin{eqnarray*}
L_{B}\left( dz^{j}\right) &=&d\bar{z}^{j} \\
TL_{B}\left( dz^{j}\right) &=&\Sigma \phi ^{jk}dz_{1}\cdots dz_{n}d\bar{z}%
_{k}.
\end{eqnarray*}
That is 
\begin{equation*}
L_{A}T\left( dz^{j}\right) =TL_{B}\left( dz^{j}\right) \text{.}
\end{equation*}
Similarly we can argue for other forms in the same way and obtain 
\begin{equation*}
L_{A}T=TL_{B}\text{.}
\end{equation*}
We can also compare $\Lambda _{A}$ and $\Lambda _{B}$ in the same way to
obtain 
\begin{equation*}
\Lambda _{A}T=T\Lambda _{B}\text{.}
\end{equation*}
Hence the variation of Hodge structure $\mathbf{sl}\left( 2\right) $ action
on $M$ was carried to the hard Lefschetz $\mathbf{sl}\left( 2\right) $
action on $W$. By symmetry, the two actions flip under the mirror
transformation $T$. $\square $

\bigskip

\bigskip

\section{\label{Automorphism}Holomorphic vs symplectic automorphisms}

\bigskip

\bigskip

\textbf{Induced holomorphic automorphisms}

For any diffeomorphism $f$ of the affine manifold $D$, its differential $df\ 
$is a diffeomorphism of $TD$ which is linear along fibers, for simplicity we
ignore the lattice $\Lambda $ in this section and write $M=TD$. We write $%
f_{B}=df:M\rightarrow M$ explicitly as 
\begin{equation*}
f_{B}\left( x^{j}+iy^{j}\right) =f^{k}\left( x^{j}\right) +i\Sigma \frac{%
\partial f^{k}}{\partial x^{j}}y^{j}.
\end{equation*}
We want to know when $f_{B}$ is a holomorphic diffeomorphism of $M$. We
compute 
\begin{eqnarray*}
&&\frac{\partial }{\partial \bar{z}^{l}}\left( f^{k}+i\Sigma \frac{\partial
f^{k}}{\partial x^{j}}y^{j}\right)  \\
&=&\frac{1}{2}\left( \frac{\partial }{\partial x^{l}}+i\frac{\partial }{%
\partial y^{l}}\right) \left( f^{k}+i\Sigma \frac{\partial f^{k}}{\partial
x^{j}}y^{j}\right)  \\
&=&\frac{1}{2}\left( \frac{\partial f^{k}}{\partial x^{l}}+i\Sigma \frac{%
\partial ^{2}f^{k}}{\partial x^{j}\partial x^{l}}y^{j}-\Sigma \frac{\partial
f^{k}}{\partial x^{j}}\delta _{l}^{j}\right)  \\
&=&\frac{i}{2}\Sigma \frac{\partial ^{2}f^{k}}{\partial x^{j}\partial x^{l}}%
y^{j}\text{.}
\end{eqnarray*}
Therefore $f_{B}$ is holomorphic on $M$ if and only if $f$ is an affine
diffeomorphism on $D$. Moreover 
\begin{equation*}
\left( f\circ g\right) _{B}=f_{B}\circ g_{B},
\end{equation*}
that is $f\rightarrow f_{B}$ is a covariant functor.

Even though $\bar{\partial}f_{B}$ does not vanish in general, its real part
does. To understand what this implies, we recall that a $T^{n}$-invariant
Calabi-Yau manifold has a natural deformation of complex structure towards
its large complex structure limit point. Its complex coordinates are given
by $dz^{j}\left( t\right) =t^{-1}dx^{j}+idy^{j}$'s with $t$ approaches $0$.
Therefore we would have 
\begin{equation*}
f_{B}\left( \frac{1}{t}x^{j}+iy^{j}\right) =\frac{1}{t}f^{k}\left(
x^{j}\right) +i\Sigma \frac{\partial f^{k}}{\partial x^{j}}y^{j},
\end{equation*}
and 
\begin{equation*}
\frac{\partial }{\partial \bar{z}^{l}\left( t\right) }\left( f^{k}+i\Sigma 
\frac{\partial f^{k}}{\partial x^{j}}y^{j}\right) =t\left( \frac{i}{2}\Sigma 
\frac{\partial ^{2}f^{k}}{\partial x^{j}\partial x^{l}}y^{j}\right) \text{.}
\end{equation*}
Namely (1) the function $f_{B}$ is holomorphic at the large complex
structure limit point $J_{\infty }$; (2) If $f$ is an affine diffeomorphism
of $D$, then $f_{B}$ is holomorphic with respect to $J_{t}$ for all $t$. We
denote this functor $f\rightarrow f_{B}$ in these two cases as follows: 
\begin{equation*}
\left( \cdot \right) _{B}:Diff\left( D\right) \rightarrow Diff\left(
M,J_{\infty }\right) ,
\end{equation*}
and 
\begin{equation*}
\left( \cdot \right) _{B}:Diff\left( D,{\scriptsize affine}\right)
\rightarrow Diff\left( M,J\right) \text{.}
\end{equation*}

\bigskip

\bigskip

\textbf{Induced symplectic automorphisms}

On the other hand, any diffeomorphism $f:D\rightarrow D$ induces a
diffeomorphism $\hat{f}:D^{\ast }\leftarrow D^{\ast }$ going the other
direction. Here $D^{\ast }\subset \mathbb{R}^{n\ast }$ denote the image of
the Legendre tranformation of $\phi $. Pulling back one forms defines a
symplectic automorphism on the total space $T^{\ast }D^{\ast }$ which is
just $M$ again (see below for explicit formula). We denote this functor as 
\begin{equation*}
\begin{array}{cccc}
\left( \cdot \right) _{A}: & Diff\left( D\right)  & \rightarrow  & 
Diff\left( M,\omega \right)  \\ 
& f & \rightarrow  & f_{A}.
\end{array}
\end{equation*}
Again 
\begin{equation*}
\left( f\circ g\right) _{A}=f_{A}\circ g_{A}.
\end{equation*}
That is $f\rightarrow f_{A}$ is also a covariant functor.

\bigskip

What if $f$ also preserves the affine structure on $D$?

The map $\hat{f}^{\ast }:M\rightarrow M$ is given by 
\begin{equation*}
\hat{f}^{\ast }\left( \hat{f}_{j}\left( x_{k}\right) ,y^{j}\right) =\left(
x_{j},\Sigma \frac{\partial \hat{f}_{k}}{\partial x_{j}}y^{k}\right) ,
\end{equation*}
for $\left( \hat{f}_{j}\left( x_{k}\right) ,y^{j}\right) \in T^{\ast
}D^{\ast }=M.$

Since $M$ is the total space of a cotangent bundle $T^{\ast }D^{\ast },$ it
has a canonical symplectic form, namely $\omega =\Sigma dx_{j}\wedge dy^{j},$
where $dx_{j}=\Sigma \phi _{jk}dx^{k}$. When the base space $D^{\ast }$ is
an affine manifold, then there is a degree two tensor $\varpi $ on its
cotangent bundle $M=T^{\ast }D^{\ast }$ whose antisymmetric part is $\omega $%
. It is given by 
\begin{equation*}
\varpi =\Sigma dx_{j}\otimes dy^{j}.
\end{equation*}
It is easy to see that $\varpi $ is well-defined on $M$.

\begin{lemma}
If $f$ is a diffeomorphism of $D$, then $f$ preserves the affine structure
on $D$ if and only if $f_{A}$ preserves $\varpi $ on $M$.
\end{lemma}

Proof: Consider the inverse of $f$ as before, $\hat{f}:D^{\ast }\rightarrow
D^{\ast }$. The pullback map it induced, $\hat{f}^{\ast }:M\rightarrow M$,
is given by 
\begin{equation*}
\hat{f}^{\ast }\left( \hat{f}_{j}\left( x_{k}\right) ,y^{j}\right) =\left(
x_{j},\Sigma \frac{\partial \hat{f}_{k}}{\partial x_{j}}y^{k}\right) ,
\end{equation*}
for $\left( \hat{f}_{j}\left( x_{k}\right) ,y^{j}\right) \in T^{\ast
}D^{\ast }=M.$ We compute 
\begin{eqnarray*}
\hat{f}^{\ast }\left( \varpi \right)  &=&\hat{f}^{\ast }\left( \Sigma
dx_{j}\otimes dy^{j}\right)  \\
&=&\Sigma dx_{j}\otimes d\left( \frac{\partial \hat{f}_{k}}{\partial x_{j}}%
y^{k}\right)  \\
&=&\Sigma dx_{j}\otimes \left( \frac{\partial \hat{f}_{k}}{\partial x_{j}}%
dy^{k}+y^{k}\frac{\partial ^{2}\hat{f}_{k}}{\partial x_{j}\partial x_{l}}%
dx_{l}\right)  \\
&=&\varpi +\Sigma y^{k}\frac{\partial ^{2}\hat{f}_{k}}{\partial
x_{j}\partial x_{l}}dx_{j}\otimes dx_{l}\text{.}
\end{eqnarray*}
Therefore $\hat{f}^{\ast }\left( \varpi \right) =\varpi $ if and only if $%
\hat{f}$ is an affine transformation on $D^{\ast }$. And this is equivalent
to $f$ being an affine transformation on $D$. $\square $

\bigskip 

We denote this functor $f\rightarrow f_{A}$ in these two cases as follows: 
\begin{equation*}
\left( \cdot \right) _{A}:Diff\left( D\right) \rightarrow Diff\left(
M,\omega \right) ,
\end{equation*}
and 
\begin{equation*}
\left( \cdot \right) _{A}:Diff\left( D,{\scriptsize affine}\right)
\rightarrow Diff\left( M,\varpi \right) \text{.}
\end{equation*}

\bigskip

\bigskip

\textbf{Transforming symplectic and holomorphic automorphisms}

We denote the spaces of those automorphisms $Diff\left( M,\ast \right) $
which are linear along fibers of the special Lagrangian fibration by $%
Diff\left( M,\ast \right) _{lin}$. Here $\ast $ may stand for $J,J_{\infty
},\omega $ or $\varpi $. Notice that for any diffeomorphism $f$ of $D$, its
induced diffeomorphisms $f_{A}$ and $f_{B}$ of $M$ are always linear along
fibers of the Lagrangian fibration $\pi :M\rightarrow D$. In fact the
converse is also true.

\begin{proposition}
(i) The map $f\rightarrow f_{B}$ induces an isomorphism, 
\begin{equation*}
\left( \cdot \right) _{B}:Diff\left( D\right) \overset{\cong }{\rightarrow }%
Diff\left( M,J_{\infty }\right) _{lin},
\end{equation*}
and similarly 
\begin{equation*}
\left( \cdot \right) _{B}:Diff\left( D,{\scriptsize affine}\right) \overset{%
\cong }{\rightarrow }Diff\left( M,J\right) _{lin}\text{.}
\end{equation*}

(ii) Moreover the map $f\rightarrow f_{A}$ induces an isomorphism, 
\begin{equation*}
\left( \cdot \right) _{A}:Diff\left( D\right) \overset{\cong }{\rightarrow }%
Diff\left( M,\omega \right) _{lin},
\end{equation*}
and similarly, 
\begin{equation*}
\left( \cdot \right) _{A}:Diff\left( D,{\scriptsize affine}\right) \overset{%
\cong }{\rightarrow }Diff\left( M,\varpi \right) _{lin}\text{.}
\end{equation*}
\end{proposition}

Proof of proposition: All these homomorphisms are obviously injective. To
prove surjectivity, we let $F$ be any diffeomorphism of $M$ which is linear
along fibers. We can write 
\begin{eqnarray*}
F &=&\left( F^{1},...,F^{n}\right) \\
F^{k} &=&f^{k}\left( x\right) +i\Sigma g_{l}^{k}\left( x\right) y^{l},
\end{eqnarray*}
for some functions $f^{k}\left( x\right) $ and $g_{l}^{k}\left( x\right) $%
's. We have 
\begin{eqnarray*}
\frac{\partial }{\partial \bar{z}^{j}}F^{k} &=&\left( \frac{\partial }{%
\partial x^{j}}+i\frac{\partial }{\partial y^{j}}\right) \left( f^{k}\left(
x\right) +i\Sigma g_{l}^{k}\left( x\right) y^{l}\right) \\
&=&\frac{\partial f^{k}}{\partial x^{j}}-g_{l}^{k}\delta _{jk}+i\frac{%
\partial g_{l}^{k}}{\partial x^{j}}y^{k}\text{.}
\end{eqnarray*}
So if $F$ preserves $J_{\infty }$, then 
\begin{equation*}
\func{Re}\frac{\partial }{\partial \bar{z}^{j}}F^{k}=0,
\end{equation*}
for all $j$ and $k$. That is $g_{j}^{k}=\frac{\partial f^{k}}{\partial x^{j}}
$, or equivalently $F=f_{B}$. If $F$ preserves $J$, additionally we have

\begin{eqnarray*}
0 &=&\frac{\partial g_{l}^{k}}{\partial x^{j}} \\
&=&\frac{\partial }{\partial x^{j}}\frac{\partial f^{k}}{\partial x^{l}}%
\text{.}
\end{eqnarray*}
That is $f$ is an affine function of $x^{j}$'s. This proves part (i). The
proof for the second part is similar. Hence we have the proposition. $%
\square $

\bigskip

Remark: Paul Yang proved that every biholomorphism of $M=TD$ with $D$ convex
is induced from an affine transformation on $D$, namely the assumption on
the function $F$ being linear along fibers is automatic. That is $Diff\left(
M,J\right) \cong Diff\left( M,J\right) _{lin}\cong Diff\left(
D,affine\right) $.

\bigskip

Next we are going to show that the mirror transformation interchanges these
two type of automorphisms. Recall that $W$ is the moduli space of special
Lagrangian tori in $M$ together with flat $U\left( 1\right) $ connections on
them. That is the moduli space of $A$-cycles $\left( C,L\right) $ with $C$ a
topological torus. Given any diffeomorphism $F:M\rightarrow M$ which is
linear along fibers of $\pi $, $F$ carries a special Lagrangian torus $C$ in 
$M$, which is a fiber to $\pi $, to another special Lagrangian torus in $M$.
The flat $U\left( 1\right) $ connection over $C$ will be carried along under 
$F$. Therefore $F$ induces a diffeomorphism of $W$, this is the mirror
transformation of $F$ and we call it $\hat{F}$ or $T\left( F\right) $.

\begin{theorem}
For $T^{n}$-invariant Calabi-Yau mirror manifolds $M$ and $W$, the above
mirror transformation $T$ induces isomorphisms : (i) 
\begin{eqnarray*}
T &:&Diff\left( M,J\right) _{lin}\overset{\cong }{\rightarrow }Diff\left(
W,\varpi \right) _{lin}, \\
T &:&Diff\left( M,\varpi \right) _{lin}\overset{\cong }{\rightarrow }%
Diff\left( W,J\right) _{lin}.
\end{eqnarray*}
and (ii) 
\begin{eqnarray*}
T &:&Diff\left( M,J_{\infty }\right) _{lin}\overset{\cong }{\rightarrow }%
Diff\left( W,\omega \right) _{lin}, \\
T &:&Diff\left( M,\omega \right) _{lin}\overset{\cong }{\rightarrow }%
Diff\left( W,J_{\infty }\right) _{lin}.
\end{eqnarray*}
Moreover the composition of two mirror transformations is the identity.
\end{theorem}

Proof of theorem: Given any $F\in Diff\left( M,J\right) _{lin}$ there is a
unique diffeomorphism $f\in Diff\left( D,affine\right) $ such that $F=f_{B}$%
. Let $\hat{f}$ be the inverse of $f$ which we consider as an affine
diffeomorphism of $D^{\ast }$. It is not difficult to verify that 
\begin{equation*}
\hat{F}=\left( \hat{f}\right) _{A}.
\end{equation*}
In particular $\hat{F}\in Diff\left( W,\omega \right) _{lin}$. Clearly 
\begin{equation*}
\widehat{\hat{F}}=F.
\end{equation*}
Other isomorphisms can be verified in the same way. Hence we have the
theorem. $\square $

\bigskip

\bigskip

\textbf{Isometries of }$M$

Recall that a diffeomorphism of a K\"{a}hler manifold preserving both the
complex structure and the symplectic structure is an isometry. Suppose $F$
is such an isometry of a $T^{n}$-invariant Calabi-Yau manifold $M,$ and we
assume that $F$ is also linear along fibers of the special Lagrangian
fibration. Then it induces a diffeomorphism $f$ of $D$ which preserves $%
g_{D}=\Sigma \phi _{jk}dx^{j}\otimes dx^{k}$. That is $f\in Diff\left(
D,g_{D}\right) $. In this case we have $F=f_{A}=f_{B}$. Hence we have 
\begin{eqnarray*}
Diff\left( D,g_{D}\right) &=&Diff\left( M,g\right) _{lin} \\
&=&Diff\left( M,J\right) _{lin}\cap Diff\left( M,\omega \right) _{lin}\text{.%
}
\end{eqnarray*}
By the above theorem, this implies that the mirror transform $\hat{F}$ lies
inside, 
\begin{equation*}
\hat{F}\in Diff\left( W,\varpi \right) _{lin}\cap Diff\left( W,J_{\infty
}\right) _{lin}\text{.}
\end{equation*}
In fact one can also show that this common intersection is simply $%
Diff\left( W,g\right) _{lin}$. A different way to see this is to observe
that the Legendre transformation from $D$ to $D^{\ast }$ preserves the
corresponding metrics $g_{D}$ and $g_{D^{\ast }}$ because 
\begin{eqnarray*}
&&\Sigma \phi ^{jk}dx_{j}\otimes dx_{k} \\
&=&\Sigma \phi ^{jk}\left( \phi _{jl}dx^{l}\right) \otimes \left( \phi
_{km}dx^{m}\right) \\
&=&\Sigma \delta _{l}^{k}\phi _{km}dx^{l}\otimes dx^{m} \\
&=&\Sigma \phi _{jk}dx^{j}\otimes dx^{k}\text{.}
\end{eqnarray*}
Therefore if $f\in Diff\left( D,g_{D}\right) $ then $\hat{f}\in Diff\left(
D^{\ast },g_{D^{\ast }}\right) $. Hence $\hat{F}$ is an isometry of $W,$ 
\begin{equation*}
\hat{F}\in Diff\left( W,g\right) .
\end{equation*}
Thus we have proved the following theorem.

\begin{theorem}
For $T^{n}$-invariant Calabi-Yau mirror manifolds $M$ and $W$, the mirror
transformation induces isomorphisms 
\begin{equation*}
T:Diff\left( M,g_{M}\right) _{lin}\overset{\cong }{\rightarrow }Diff\left(
W,g_{W}\right) _{lin}\text{.}
\end{equation*}
Moreover the composition of two mirror transformations is the identity.
\end{theorem}

\bigskip

\section{A- and B-connections}

Mirror transformations of other A- and B-cycles can be interpreted as
generalization of the classical duality between Blaschke connection and its
conjugate connection via Legendre transformation. We first recall these
classical geometries.

\bigskip

\textbf{Blaschke connection and its conjugate connection}

On a special affine manifold $D$, $\phi $ is a section of a trivial real
line bundle over $D$. For simplicity we assume $D\subset \mathbb{R}^{n}$ and 
$\phi $ is a convex function on $D$. If $G\subset D\times \mathbb{R}$ denote
the graph of $\phi $. Using the affine structure on $D\times \mathbb{R}$ one
can define an affine normal $\nu $ which is a transversal vector field along 
$G$ (see for example \cite{CY}): If we parallel translate the tangent plane
of $G$, its intersection with $G$ determines a small convex domain. Its
center of gravity then traces out a curve in the space whose initial
direction is the affine normal direction.

Using $\nu $, we can decompose the restriction of the standard affine
connection on $\mathbb{R}^{n}\times \mathbb{R}$ to $G$ into tangent
directions and normal direction. So we obtain an induced torsion free
connection on $G$, called the Blaschke connection \cite{Bl}, or the
B-connection, $\!_{B}\!\nabla $. The convexity of $\phi $ implies that the
second fundamental form is a positive definite symmetric two tensor $g_{G}$
on $G$. Its Levi-Civita connection is denoted as $\nabla ^{LC}$. We define a
conjugate connection $\!_{A}\!\nabla $ by 
\begin{equation*}
Xg_{G}\left( Y,Z\right) =g_{G}\left( \!_{B}\!\nabla _{X}Y,Z\right)
+g_{G}\left( Y,\!\,_{A}\!\nabla _{X}Z\right) \text{.}
\end{equation*}
We call it an A-connection. The two connections $_{A}\!\nabla $ and $%
\!_{B}\!\nabla $ on $D$ induce torsion free connections on $M=TD$ by
pullback. We continue to call them A-connection $\!_{A}\!\nabla $ and
B-connection $\!_{B}\!\nabla $.

\bigskip

When the function $\phi $ satisfying the real Monge-Amper\'{e} equation,
then $G$ is a parabolic affine sphere in $\mathbb{R}^{n}\times \mathbb{R}$.
Namely the affine normal $\nu $ of $G$ in $\mathbb{R}^{n}\times \mathbb{R}$
is the unit vector along the last direction, that is the fiber direction of
the real line bundle over $D$. By abuse of notations, we identify $G$ with $%
D $ via the projection to the first factor in $D\times \mathbb{R}$. These
two torsion free connections $_{A}\,\nabla $ and $\,\!_{B}\!\nabla $ on $D$
are flat in this case. In term of the affine coordinates $x^{j}$'s on $D$, $%
\,$the B-connection $\!_{B}\!\nabla $ is just given by the exterior
differentiation $d$. The A-connection is 
\begin{equation*}
\,\!_{A}\!\nabla =d+\Sigma \phi ^{jk}\phi _{klm}\text{.}
\end{equation*}
One can check directly that it has zero curvature. We will see later that
this also follows from the Legendre transformation (or the mirror
transformation).

\bigskip

\bigskip

\textbf{A- and B-connections on $T^{n}$-invariant Calabi-Yau manifolds}

From above, we have two torsion free flat connections $_{A}\,\nabla $ and $%
\,\!_{B}\!\nabla $ on $M=TD/\Lambda $. Recall that the complex structure on $%
M$ is given $z^{j}=x^{j}+iy^{j}$ and its symplectic form is $\omega
_{M}=\Sigma \phi _{jk}\left( x\right) dx^{j}\wedge dx^{k}$. Since $%
\,\!_{B}\!\nabla $ is the same as the exterior differentiation on the affine
coordinates $x^{j}$'s on $D$, it preserves the complex structure on $M$. In
fact $_{A}\nabla $ preserves the symplectic structure on $M$.

\begin{proposition}
Let $M=TD/\Lambda $ be a $T^{n}$-invariant Calabi-Yau manifold as before.
Then its A-connection $_{A}\,\nabla $ and B-connection $\!_{B}\!\nabla $
satisfies 
\begin{eqnarray*}
\medskip _{A}\nabla \omega _{M} &=&0, \\
\!_{B}\!\nabla J &=&0.
\end{eqnarray*}
\end{proposition}

Proof of proposition: We have seen that $\!_{B}\!\nabla J=0$. From previous
discussions, $_{A}\triangledown $ is a torsion free flat connection on $M$.
To check that it preserves the symplectic form, we recall that $\omega
_{M}=\Sigma \phi _{jk}\left( x\right) dx^{j}\wedge dy^{k}$ and $%
_{A}\!\triangledown =d+\Gamma _{lm}^{j}dx^{m}$ where $\Gamma _{lm}^{j}=\phi
^{jk}\phi _{klm}$. Note that $_{A}\!\triangledown \left( dy^{k}\right) =0$
because $_{A}\!\triangledown $ is induced from the base $D$. Therefore 
\begin{eqnarray*}
_{A}\triangledown \omega _{M} &=&\Sigma \phi _{jkl}dx^{l}\otimes \left(
dx^{j}\wedge dy^{k}\right) -\Sigma \phi _{jk}\Gamma _{pq}^{j}dx^{q}\otimes
\left( dx^{p}\wedge dy^{k}\right) \\
&=&\Sigma \phi _{jkl}dx^{l}\otimes \left( dx^{j}\wedge dy^{k}\right) -\Sigma
\phi _{jk}\phi ^{jl}\phi _{lpq}dx^{q}\otimes \left( dx^{p}\wedge
dy^{k}\right) \\
&=&\Sigma \phi _{jkl}dx^{l}\otimes \left( dx^{j}\wedge dy^{k}\right) -\Sigma
\phi _{kpq}dx^{q}\otimes \left( dx^{p}\wedge dy^{k}\right) \\
&=&0.
\end{eqnarray*}
We have use the symmetry of $\phi _{jkl}$ with respect to its indices. $%
\square $

\bigskip

In section \ref{TransformCycle}, we will see that the mirror transformation
of flat connection $_{A}\triangledown $ (resp. $_{B}\triangledown $) on $M$
is the flat connection $_{B}\triangledown $ (resp. $_{A}\triangledown $) on
the zero section in $W$ and vice versa. In particular the Levi-Civita
connection $\triangledown ^{LC}=\left( _{A}\triangledown
+\,_{B}\triangledown \right) /2$ is preserved under the mirror
transformation. In fact this is a special case of the mirror transformation
between A- and B-cycles on $M$ and $W$.

\bigskip

\section{\label{TransformCycle}Transformation of A- and B-cycles}

In this section we discuss how certain A-cycles on $M$ will transform to
B-cycles on $W$. This materials is largely borrowed from \cite{LYZ}.\ 

\bigskip

\textbf{Transforming A- and B-connections}

Recall the A-connection $_{A}\triangledown $ on $M$ is $d+\Sigma \Gamma
_{kl}^{j}dx^{l}$, where $\Gamma _{kl}^{j}=\Sigma \phi ^{jm}\phi _{mkl}$, in
the affine coordinate system. Let us consider $M$ and $W$ with their dual
special Lagrangian tori fibrations. The restriction of $_{A}\triangledown $
on each fiber in $M$ is trivial because $dx^{l}$'s vanish along fiber
directions. Since the dual torus $T^{\ast }$ parametrizes flat $U\left(
1\right) $ connections on $T$, the restriction of $_{A}\triangledown $
corresponds to the origin of the corresponding dual torus. Putting all
fibers on $M$ together, we obtain the zero section in $W$. This is the
Fourier transformation.

However this is not the end of the story, the second fundamental form of $%
_{A}\triangledown $ on each fiber in $M$ is non-trivial. This induces a
connection on the zero section in $W$. To determine this connection, we need
to perform the Legendre transformation on $M$, 
\begin{eqnarray*}
_{A}\nabla _{\frac{\partial }{\partial x^{j}}}\left( \tfrac{\partial }{%
\partial x^{k}}\right) &=&\Gamma _{jk}^{l}\frac{\partial }{\partial x^{l}} \\
\nabla _{\frac{\partial }{\partial x^{j}}}\left( \Sigma \phi _{kq}\tfrac{%
\partial }{\partial x^{q}}\right) &=&\Sigma \phi ^{lm}\phi _{mjk}\frac{%
\partial }{\partial x^{l}} \\
\Sigma \phi _{kqj}\tfrac{\partial }{\partial x_{q}}+\Sigma \phi _{jp}\phi
_{kq}\nabla _{\frac{\partial }{\partial x_{p}}}\left( \tfrac{\partial }{%
\partial x_{q}}\right) &=&\Sigma \phi _{mjk}\frac{\partial }{\partial x_{m}}
\\
\Sigma \phi _{jp}\phi _{kq}\nabla _{\frac{\partial }{\partial x_{p}}}\left( 
\tfrac{\partial }{\partial x_{q}}\right) &=&0.
\end{eqnarray*}
That is 
\begin{equation*}
\triangledown _{\frac{\partial }{\partial x_{p}}}\left( \tfrac{\partial }{%
\partial x_{q}}\right) =0,
\end{equation*}
or equivalently the induced connection on the zero section of $W$ is $d$ in
the affine coordinate system of $W$. This is exactly the B-connection $%
_{B}\triangledown $.

\bigskip

Conversely if we start with the B-connection $_{B}\triangledown $ on the
whole manifold $M$, its mirror transformation will be the A-connection $%
_{A}\triangledown $ on the zero section of $W$. In particular we recover the
classical duality between the Blaschke connection and its conjugate
connection for the parabolic affine sphere. Such duality is in fact more
interesting for other affine hypersurfaces (see for example \cite{Lo}).
Summarizing we have the following theorem.

\begin{theorem}
For a $T^{n}$-invariant manifold $M$, the above mirror transformation take
the A-connection (resp. B-connection) on the whole space $M$ to the
B-connection (resp. A-connection) on the zero section of $W$.
\end{theorem}

\bigskip

\textbf{Transforming special Lagrangian sections}

Now we are going to generalize the previous picture to duality between other
supersymmetric cycles. Let $\left( C,E\right) $ be an A-cycle in $M$ such
that $C$ is a section of the special Lagrangian fibration $\pi :M\rightarrow
D$.

Note that $C=\left\{ y=y\left( x\right) \right\} \subset M$ being Lagrangian
with respect to $\omega _{M}=\Sigma \phi _{jk}dx^{j}dy^{k}$ is equivalent to

\begin{equation*}
\frac{\partial }{\partial x^{j}}(y^{l}\phi _{lk})=\frac{\partial }{\partial
x^{k}}(\phi _{lj}y^{l}).
\end{equation*}
Therefore locally there is a function $f$ on $D$ such that 
\begin{equation*}
y^{j}=\Sigma \phi ^{jk}\frac{\partial f}{\partial x^{k}}.
\end{equation*}
Next we want to understand the special condition on $C$. Namely 
\begin{equation*}
\func{Im}e^{i\theta }\Omega _{M}|_{C}=0.
\end{equation*}
Recall that the holomorphic volume form on $M$ equals $\Omega
_{M}=dz^{1}\wedge dz^{2}\wedge ...\wedge dz^{n}$.

On the Lagrangian section $C$ we have 
\begin{eqnarray*}
dy^{j} &=&d\left( \Sigma \phi ^{jk}\frac{\partial f}{\partial x^{k}}\right) 
\\
&=&\Sigma \phi ^{jl}\left( \frac{\partial ^{2}f}{\partial x^{l}\partial x^{k}%
}-\phi ^{pq}\phi _{lkp}\frac{\partial f}{\partial x^{q}}\right) dx^{k} \\
&=&\Sigma \phi ^{jl}\,_{A}Hess\left( f\right) _{lk}dx^{k}\text{.}
\end{eqnarray*}
Here $_{A}Hess\left( f\right) $ denote the Hessian of $f$ with respect to
the restriction of the torsion free A-connection $_{A}\triangledown $ and we
use the affine coordinate on $D$ to parametrize the section $C$. We have 
\begin{eqnarray*}
dz^{j} &=&dx^{j}+idy^{j} \\
&=&\Sigma \left( \delta _{jk}+i\phi ^{jl}\left( \frac{\partial ^{2}f}{%
\partial x^{l}\partial x^{k}}-\phi ^{pq}\phi _{lkp}\frac{\partial f}{%
\partial x^{q}}\right) \right) dx^{k},
\end{eqnarray*}
and 
\begin{eqnarray*}
\Omega _{M}|_{C} &=&\det \left( I+ig^{-1}\,_{A}Hess\left( f\right) \right)
dx^{1}\wedge ...\wedge dx^{n} \\
&=&\det \left( g\right) ^{-1}\det \left( g+i\,_{A}Hess\left( f\right)
\right) dx^{1}\wedge ...\wedge dx^{n},
\end{eqnarray*}
Hence $C$ is a special Lagrangian section if and only if 
\begin{equation*}
\func{Im}e^{i\theta }\det \left( g+i\,_{A}Hess\left( f\right) \right) =0.
\end{equation*}

{}Now we perform the fiberwise Fourier transformation on $M$. On each torus
fiber $T$, the special Lagrangian section $C$ determines a point $y=\left(
y^{1},...,y^{n}\right) $ on it, and therefore a flat $U\left( 1\right) $
connection $D_{y}$ on its dual torus $T^{\ast }$. Explicitly, we have 
\begin{equation*}
D_{y}=d+i\Sigma y^{j}dy_{j}.
\end{equation*}

By putting all these fibers together, we obtain a $U\left( 1\right) $
connection $\triangledown _{A}$ on the whole $W$, 
\begin{equation*}
\triangledown _{A}=d+i\Sigma y^{j}dy_{j}.
\end{equation*}
Its curvature two form is given by, 
\begin{equation*}
F_{A}=\left( \triangledown _{A}\right) ^{2}=\Sigma i\frac{\partial y^{j}}{%
\partial x_{k}}dx_{k}\wedge dy_{j}.
\end{equation*}
The $\left( 2,0\right) $ component of the curvature equals 
\begin{equation*}
F_{A}^{2,0}=\frac{1}{2}\Sigma \left( \frac{\partial y^{k}}{\partial x_{j}}-%
\frac{\partial y^{j}}{\partial x_{k}}\right) dz_{j}\wedge dz_{k}.
\end{equation*}
Therefore $\triangledown _{A}$ gives a holomorphic line bundle on $W$ if and
only if 
\begin{equation*}
\frac{\partial y^{k}}{\partial x_{j}}=\frac{\partial y^{j}}{\partial x_{k}},
\end{equation*}
for all $j,k$. This is equivalent to the existence of a function $f=f\left(
x_{j}\right) $ on $D$ such that 
\begin{equation*}
y^{j}=\frac{\partial f}{\partial x_{j}}.
\end{equation*}
Therefore we can rewrite the curvature tensor as 
\begin{equation*}
F_{A}=i\Sigma \,_{B}Hess\left( f\right) _{jk}dx_{k}\wedge dy_{j}\text{.}
\end{equation*}
Here $_{B}Hess\left( f\right) $ is the Hessian of $f$ with respect to the
B-connection $_{B}\triangledown $ on $W$.

To compare with the $M$ side, we use the Legendre transformation to write 
\begin{equation*}
y^{j}=\Sigma \phi ^{jk}\frac{\partial f}{\partial x^{k}}.
\end{equation*}
Then $_{B}Hess\left( f\right) $ on $W$ becomes $_{A}Hess\left( f\right) $ on 
$M.$ Therefore the cycle $C\subset M$ being a special Lagrangian is
equivalent to 
\begin{eqnarray*}
F_{A}^{2,0} &=&0, \\
\func{Im}e^{i\theta }\left( \omega _{W}+F_{A}\right) ^{n} &=&0.
\end{eqnarray*}

Next we bring back the flat $U\left( 1\right) $ connection on $E$ over $C$
to the picture. We still use the affine coordinates on $D$ to parametrize $C$
because it is a section. We can express the flat connection on $C$ as

\begin{equation*}
d+ide=d+i\Sigma \frac{\partial e}{\partial x^{k}}dx^{k}
\end{equation*}
for some function $e=e\left( x\right) $ on $C$. Now this connection will be
added to the previous one on $W$ as the second fundamental form along
fibers. We still call this connection $\triangledown _{A}.$ We have 
\begin{eqnarray*}
\triangledown _{A} &=&d+i\Sigma y^{j}dy_{j}+ide \\
&=&d+i\Sigma \phi ^{jk}\frac{\partial f}{\partial x^{k}}dy_{j}+i\Sigma \frac{%
\partial e}{\partial x_{j}}dx_{j}.
\end{eqnarray*}
It is easy to see that the curvature form of this new connection is the same
as the old one. In particular the transformed connection $\triangledown _{A}$
on $W$ continues to satisfy the deformed Hermitian-Yang-Mills equations. 
\begin{eqnarray*}
F_{A}^{0,2} &=&0, \\
\func{Im}e^{i\theta }\left( \omega _{W}+F_{A}\right) ^{m} &=&0,
\end{eqnarray*}
Therefore the mirror transformation of the A-cycle $\left( C,E\right) $ on $%
M $ produces a B-cycle on $W$. The same approach work for higher rank
unitary bundle over the section $C$. This transformation is explained with
more details in \cite{LYZ}.

\bigskip

\bigskip

\textbf{Transforming graded tangent spaces}

Recall from \cite{L1} that the tangent space of the moduli space of A-cycle $%
\left( C,E\right) $ in $M$ is the space of complex harmonic one form with
valued in the adjoint bundle. That is 
\begin{equation*}
T\left( _{A}\mathcal{M}\left( M\right) \right) =H^{1}\left( C,ad\left(
E\right) \right) \otimes \mathbb{C}\text{.}
\end{equation*}
And the tangent space of the moduli space of B-cycle $\left( C,E\right)
=\left( W,E\right) $ in $W$ is the space of deformed $\bar{\partial}$%
-harmonic one form with valued in the adjoint bundle. 
\begin{equation*}
T\left( _{B}\mathcal{M}\left( W\right) \right) =QH^{1}\left( C,End\left(
E\right) \right) \text{.}
\end{equation*}

A form $B\in \Omega ^{0,q}\left( C,End\left( E\right) \right) $ is called a
deformed $\bar{\partial}$-harmonic form if it satisfies the following
deformation of the harmonic form equations: 
\begin{eqnarray*}
\bar{\partial}B &=&0, \\
\func{Im}e^{i\theta }\left( \omega +F\right) ^{m-q}\wedge \partial B &=&0.
\end{eqnarray*}
Here $m$ is the complex dimension of $C$.

The graded tangent spaces are given by 
\begin{eqnarray*}
T^{graded}\left( _{A}\mathcal{M}\left( M\right) \right)  &=&\oplus
_{k}H^{k}\left( C,ad\left( E\right) \right) \otimes \mathbb{C}, \\
T^{graded}\left( _{B}\mathcal{M}\left( W\right) \right)  &=&\oplus
_{k}QH^{k}\left( C,End\left( E\right) \right) \text{.}
\end{eqnarray*}

Now we identify these two spaces when $C\subset M$ is a special Lagrangian
section. It is easy to see that the linearization of the above
transformation of A-cycles on $M$ to B-cycles on $W$ is the following
homomorphism

\begin{eqnarray*}
\Omega ^{1}\left( C,ad\left( E\right) \right) \otimes \mathbb{C}
&\rightarrow &\Omega ^{0,1}\left( W,End\left( E\right) \right) \\
dx^{j} &\rightarrow &\Sigma \frac{i}{2}\phi ^{jk}d\bar{z}_{k}.
\end{eqnarray*}
We extend that homomorphism to higher degree forms, in the obvious way, 
\begin{equation*}
\Omega ^{q}\left( C,ad\left( E\right) \right) \otimes \mathbb{C}\rightarrow
\Omega ^{0,q}\left( W,End\left( E\right) \right) .
\end{equation*}

It is verified in \cite{LYZ} that the harmonic form equation on $\Omega
^{q}\left( C,ad\left( E\right) \right) \otimes \mathbb{C}$ is transformed to
the deformed harmonic form equation on $\Omega ^{0,q}\left( W,End\left(
E\right) \right) .$ Namely the image of $H^{q}\left( C,ad\left( E\right)
\right) \otimes \mathbb{C}$ under the above homomorphism is inside $%
QH^{q}\left( W,End\left( E\right) \right) $. In fact the image is given
precisely by those forms which are invariant along fiber directions.

As a corollary of this identification, we can also see that the mirror
transformation between moduli space of cycles, $_{A}\mathcal{M}\left(
M\right) \rightarrow _{B}\mathcal{M}\left( W\right) $, is a holomorphic map.

\bigskip

\textbf{Identifying correlation functions}

The correlation functions on these moduli spaces of cycles are certain
n-forms on them (see for example \cite{L1} for the intrinsic definition). On
the $M$ side, it is given by 
\begin{equation*}
_{A}\Omega \left( C,E\right) \left( \alpha _{1},...,\alpha _{n}\right)
=\int_{C}Tr_{E}\left[ \alpha _{1}\wedge ...\wedge \alpha _{n}\right] _{sym},
\end{equation*}
for $\alpha _{j}\in \Omega ^{1}\left( C,ad\left( E\right) \right) \otimes 
\mathbb{C}$ at a A-cycle $\left( C,E\right) $. On the $W$ side, it is given
by

\begin{equation*}
_{B}\Omega \left( C,E\right) \left( \beta _{1},...,\beta _{n}\right)
=\int_{W}\Omega _{W}Tr_{E}\left[ \beta _{1}\wedge \cdots \wedge \beta _{n}%
\right] _{sym},
\end{equation*}
for $\beta _{j}\in \Omega ^{0,1}\left( W,End\left( E\right) \right) $ at a
B-cycle $\left( C,E\right) =\left( W,E\right) $. If $C\neq W,$ then the
formula is more complicated (see \cite{L1}).

One can verify directly that the n-form $_{B}\Omega $ on the $W$ side is
pullback to $_{A}\Omega $ on the $M$ side under the above mirror
transformation (see \cite{LYZ} for details). This verifies Vafa conjecture
for rank one bundles in the $T^{n}$-invariant Calabi-Yau case. His
conjecture says that the moduli spaces of A- and B-cycles, together with
their correlation functions, on mirror manifolds should be identified. In
general this identification should require instanton corrections.

\section{$T^{n}$-invariant hyperk\"{a}hler manifolds}

A Riemannian manifold $M$ of dimension $4n$ with holonomy group equals $%
Sp\left( n\right) \subset SU\left( 2n\right) $ is called a hyperk\"{a}hler
manifold.

\bigskip

\textbf{$T^{n}$-invariant hyperk\"{a}hler manifolds}

As we discussed in the $T^{n}$-invariant Calabi-Yau manifolds, let $D$ be an
affine manifold with local coordinates $x^{j}$'s and $\phi \left( x\right) $
be a solution to the real Monge-Amper\'{e} equation $\det \left( \frac{%
\partial ^{2}\phi }{\partial x^{j}\partial x^{k}}\right) =1$. Then both its
tangent bundle $TD$ and cotangent bundle $T^{\ast }D$ are naturally $T^{n}$%
-invariant Calabi-Yau manifolds. Moreover they are mirror to each other. If
we denote the local coordinate of $TD$ as $x^{j}$ and $y^{j}$'s. Then the
complex structure of $TD$ is determined by $dx^{j}+idy^{j}$'s as being $%
\left( 1,0\right) $ forms and we call this complex structure $J$. Its
symplectic form is given by $\omega =\Sigma \phi _{jk}\left( x\right)
dx^{j}\wedge dy^{k}$ and its Ricci flat metric is $g=\Sigma \phi _{jk}\left(
dx^{j}\otimes dx^{k}+dy^{j}\otimes dy^{k}\right) .$

Now we consider its cotangent bundle $M=T^{\ast }\left( TD\right) $ and
denote the dual coordinates for $x^{j}$ and $y^{j}$ as $u_{j}$ and $v_{j}$
respectively. Therefore the induced metric on $M$ is given by 
\begin{equation*}
g_{M}=\Sigma \phi _{jk}\left( dx^{j}\otimes dx^{k}+dy^{j}\otimes
dy^{k}\right) +\Sigma \phi ^{jk}\left( du_{j}\otimes du_{k}+dv_{j}\otimes
dv_{k}\right) \text{,}
\end{equation*}
and its induced complex structure $J$ is determined by $dx^{j}+idy^{j}$'s
and $du_{j}-idv_{j}$'s as being $\left( 1,0\right) $ forms. Its
corresponding symplectic form $\omega _{J}$ is given by 
\begin{equation*}
\omega _{J}=\Sigma \phi _{jk}dx^{j}\wedge dy^{k}-\Sigma \phi
^{jk}du_{j}\wedge dv_{k}.
\end{equation*}

Since $M$ is the cotangent bundle of a complex manifold, it has a natural
holomorphic symplectic form which we denote as $\eta _{J}$ and it is given
by 
\begin{equation*}
\eta _{J}=\Sigma \left( dx^{j}+idy^{j}\right) \wedge \left(
du_{j}-idv_{j}\right) \text{.}
\end{equation*}
Notice that the projection $\pi :M\rightarrow TD$ is a holomorphic
Lagrangian fibration with respect to $\eta _{J}$.

We are going see that $M$ carries a natural hyperk\"{a}hler structure. If we
denote the real and imaginary part of $\eta _{J}$ by $\omega _{I}$ and $%
\omega _{K}$ respectively, then they are both real symplectic form on $M$.
Explicitly we have 
\begin{eqnarray*}
\omega _{I} &=&\func{Re}\eta _{J}=\Sigma \left( dx^{j}\wedge
du_{j}+dy^{j}\wedge dv_{j}\right) , \\
\omega _{K} &=&\func{Im}\eta _{J}=\Sigma \left( dx^{j}\wedge
dv_{j}-dy^{j}\wedge du_{j}\right) \text{.}
\end{eqnarray*}
They determine almost complex structures $I$ and $K$ on $M$ respectively. In
fact these are both integrable complex structures. If we use the following
change of variables, $du^{j}=\phi ^{jk}du_{k}$ and $dv^{j}=\phi ^{jk}dv_{k}$
then the complex structure of $I$ is determined by $dx^{j}+idu^{j}$ and $%
dy^{j}+idv^{j}$ as being $\left( 1,0\right) $ forms. Similarly the complex
structure of $K$ is determined by $dx^{j}+idv^{j}$ and $dy^{j}-idu^{j}$ as
being $\left( 1,0\right) $ forms. It follows from direct calculations that
both $\left( M,g,I,\omega _{I}\right) $ and $\left( M,g,K,\omega _{K}\right) 
$ are Calabi-Yau structures on $M$. We can easily verify the following lemma.

\begin{lemma}
$I^{2}=J^{2}=K^{2}=IJK=-id$. Namely $\left( M,g\right) $ is a hyperk\"{a}%
hler manifold.
\end{lemma}

Remark: We call such $M$ a $T^{n}$-invariant hyperkahler manifold. Instead
of $T^{\ast }\left( TD\right) $ we can also consider $T\left( T^{\ast
}D\right) $ and it also has a natural hyperk\"{a}hler structure constructed
in a similar way. In fact these two are isomorphic hyperk\"{a}hler manifolds.

\bigskip

\textbf{An }$\mathbf{so}\left( 4,1\right) $\textbf{\ action on cohomology}

For a hyperk\"{a}hler manifold $M$, there is a $S^{2}$-family of K\"{a}hler
structures $\omega _{t}$ on it: For any $t=\left( a,b,c\right) \in \mathbb{R}%
^{3}$ with $a^{2}+b^{2}+c^{2}=1$, $\omega _{t}=a\omega _{I}+b\omega
_{J}+c\omega _{K}$ is a K\"{a}hler metric on $M$. For each $\omega _{t}$,
there is a corresponding hard Lefschetz $\mathbf{sl}\left( 2\right) $ action
on its cohomology group $H^{\ast }\left( M,\mathbb{R}\right) $. It is showed
by Verbitsky in \cite{Ve} that this $S^{2}$ family of $\mathbf{sl}\left(
2\right) $ actions on $H^{\ast }\left( M,\mathbb{R}\right) $ in fact
determines an $\mathbf{so}\left( 4,1\right) $ action on cohomology. It is
interesting to compare the $\mathbf{so}\left( 3,1\right) $ action from
Gopakumar-Vafa conjecture with this $\mathbf{so}\left( 4,1\right) $ action
when $M$ admits a holomorphic Lagrangian fibration.

Note that $\mathbf{sl}\left( 2\right) \mathbf{=so}\left( 2,1\right) $\textbf{%
\ }and\textbf{\ }$\mathbf{sl}\left( 2\right) \mathbf{\times sl}\left(
2\right) \mathbf{=so}\left( 3,1\right) .$ Therefore the cohomology group of
K\"{a}hler manifolds admit $\mathbf{so}\left( 2,1\right) $ actions, the
cohomology of $T^{n}$-invariant Calabi-Yau manifolds admit $\mathbf{so}%
\left( 3,1\right) $ actions and the cohomology of hyperk\"{a}hler manifolds
admit $\mathbf{so}\left( 4,1\right) $ actions. We are going to show that the 
$\mathbf{so}\left( 3,1\right) $ action we constructed in the $T^{n}$%
-invariant Calabi-Yau case is naturally embedded inside this $\mathbf{so}%
\left( 4,1\right) $ action for hyperk\"{a}hler manifolds. This is analogous
to the statement that the hard Lefschetz $\mathbf{so}\left( 2,1\right) $
action for K\"{a}hler manifolds is part of the $\mathbf{so}\left( 3,1\right) 
$ action for Calabi-Yau manifolds, at least in the $T^{n}$-invariant case.

\bigskip

\textbf{Embedding }$\mathbf{sl}\left( 2\right) \times \mathbf{sl}\left(
2\right) $\textbf{\ inside hyperk\"{a}hler }$\mathbf{so}\left( 4,1\right) $ 
\textbf{action}

As we discussed before, besides the hard Lefschetz $\mathbf{sl}\left(
2\right) $ action on $\Omega ^{\ast ,\ast }\left( M\right) $, the other $%
\mathbf{sl}\left( 2\right) $ action comes from a variation of complex
structure on $M$. For our $T^{n}$-invariant hyperk\"{a}hler manifold $M$ as
above with the complex and K\"{a}hler structure $I$ and $\omega _{I}$ and
special Lagrangian fibration $\pi :M\rightarrow TD$, the second $\mathbf{sl}%
\left( 2\right) $ action on $M$ can be expressed using $\omega _{J}$ and $%
\omega _{K}$. That is we have a natural embedding of the $\mathbf{sl}\left(
2\right) \mathbf{\times sl}\left( 2\right) $ action into the hyperk\"{a}hler 
$\mathbf{so}\left( 4,1\right) $ action on $M$.

To verify this, we recall that the operator $L_{B}$ in the $\mathbf{sl}%
\left( 2\right) $ action coming from the VHS will send $dx^{j}+idu^{j}$ (we
write $du^{j}=\phi ^{jk}du_{k}$) to $dx^{j}-idu^{j}$. On the other hand, for
the operators $L_{J},\Lambda _{K}$ in the hyperk\"{a}hler\textbf{\ }$\mathbf{%
so}\left( 4,1\right) $ action, we have 
\begin{eqnarray*}
&&\left[ L_{J},\Lambda _{K}\right] \left( dx^{j}+idu^{j}\right)  \\
&=&-\Lambda _{K}L_{J}\left( dx^{j}+idu^{j}\right)  \\
&=&-\Lambda _{K}\left( dx^{j}+idu^{j}\right) \left( \Sigma \phi
_{kl}dx^{k}dy^{l}-\phi ^{kl}du_{k}dv_{l}\right)  \\
&=&-\left( \phi ^{kj}du_{k}-idx^{k}\right)  \\
&=&i\left( dx^{j}-idu^{j}\right) ,
\end{eqnarray*}
because $\omega _{K}=\Sigma \left( dx^{j}dv_{j}-dy^{j}du_{j}\right) $. The
same holds true for all other forms. Thus we have the following theorem.

\begin{theorem}
For any $T^{n}$-invariant hyperk\"{a}hler manifold $M$, its Calabi-Yau $%
\mathbf{so}\left( 3,1\right) =\mathbf{sl}\left( 2\right) \times \mathbf{sl}%
\left( 2\right) $ action on cohomology embeds naturally inside the hyperk%
\"{a}hler $\mathbf{so}\left( 4,1\right) $ action.
\end{theorem}

\end{document}